\documentclass[oneside,reqno]{amsart}
\usepackage{amsmath}
\usepackage{pspicture} 
\usepackage{verbatim}
\usepackage{eucal} 
\usepackage{latexsym}
\usepackage{amsfonts,amssymb,amsxtra}

\usepackage{amscd}
\usepackage{amscd}

\usepackage{fancybox,ascmac}

\usepackage{latexsym}

\newtheorem{theorem}{Theorem}[section]

\newtheorem{claim}[theorem]{Claim}

\newtheorem{corollary}[theorem]{Corollary}

\newtheorem{lemma}[theorem]{Lemma}

\newtheorem{proposition}[theorem]{Proposition}

\theoremstyle{definition}%Theorems et al in roman font.

\newtheorem{definition}[theorem]{Definition}
\newtheorem{example}[theorem]{Example}

\newtheorem{remark}[theorem]{Remark}

\newcommand\supp{\operatorname{supp}}

\newcommand\red{\text{red\ \!}}

\newcommand\spa{\text{span\ \!}}

\newcommand\ad{\text{ad\ \!}}

\newcommand\sll{\operatorname{sl}}

\begin{document}

%\today

\enskip

\title{New Lie tori from Naoi tori}

\author{Yoji Yoshii}

\dedicatory{Dedicated to Professor Jun Morita
on the occasion of his 60th birthday
}

\address[Yoji Yoshii]
{Iwate University\\
3-18-33 Ueda, Morioka, Iwate, Japan 020-8550}
\email{yoshii@iwate-u.ac.jp}
%17A30 Algebras satisfying other identities
%17A75 Composition algebras
%17B60 Lie (super)algebras associated with other structures (associative, Jordan, etc.)
%17C40 Exceptional Jordan structures
%17B67 Kac-Moody (super)algebras (structure and representation theory)
\date{}

\begin{abstract}
We define general Lie tori
which generalize original Lie tori.
We show that a Naoi torus is a general Lie torus.
We give examples and prove several properties of general Lie tori.
We also review isotopies of Lie tori, and prove that a general Lie torus is, in fact, isotopic to an original Lie torus.
Finally, we suggest a very simple way of defining a Lie torus 
corresponding to a locally extended affine root system $\frak R$,
which we call a Lie $\frak R$-torus.
\end{abstract}

\maketitle

Throughout the paper $F$ is a field of characteristic 0.
For a subset $S$ of an abelian group,
the subgroup generated by $S$ is denoted by
$\langle S\rangle$.

\enskip

\enskip

\section{Introduction}
\label{sec:intro}

Naoi showed  in \cite{Na} that Lie tori are not enough to describe the fixed algebras of
multi-loop algebras studied in \cite{ABFP}.
Thus he defined a modified Lie torus in \cite{Na}, which we call a {\bf Naoi torus}
(see Definition \ref{defNaoi}).
We define a new wider class of Lie tori, called
{\bf general Lie tori} 
by a slight modification of the definition of original Lie tori
(see Definition \ref{defgeneral}).
We show that any Naoi torus is a general Lie torus in Theorem \ref{naoiisgen}.

Next we review the notion of isotopies of Lie tori,
introduced in \cite{AF}.
Let us call an original Lie torus defined in \cite{Ne} or \cite{Y2}
a {\bf normal Lie torus}.
As very simple examples,
the loop algebra $\sll_2(F[t^{\pm 1}])$ is a normal Lie torus,
and the subalgebra $P$ of $\sll_2(F[t^{\pm 1}])$
generated by $e\otimes t$, $f\otimes t^{-1}$ and $h\otimes t^{\pm 2}$,
where
$e=
\begin{pmatrix}
0&1\\
0&0
\end{pmatrix}
$,
$f=
\begin{pmatrix}
0&0\\
1&0
\end{pmatrix}
$
and
$h=
\begin{pmatrix}
1&0\\
0&-1
\end{pmatrix}
$,
is not normal but a general Lie torus.
We will show that
any general Lie torus is isotopic to a normal Lie torus.
Thus some important properties of normal Lie tori still hold for general Lie tori.
For example, 
we show that there exists a nonzero symmetric invariant graded form
on a general Lie $n$-torus in Corollary \ref{efrom}.

However,
the support for the grading of a general Lie torus is quite different from the normal case.
We study each support, which is an example of a so-called  {\bf reflection space} $E$ of an abelian group $G$,
i.e.,
a subset $E$ of $G$ satisfies $2x-y\in E$ for all $x,y\in E$.
A reflection space $S$ of the normal case 
has a stronger condition,
namely,
$S$ satisfies
$x-2y\in S$ for all $x,y\in S$,
which is called a {\bf symmetric reflection space},
as in \cite{LN2}.
A reflection space and a symmetric reflection space look similar,
but they are very different.
For example, $m\Bbb Z+n$ for any $m, n\in\Bbb Z$ is
a reflection space of $\Bbb Z$,
but it is symmetric only when $n=0$,
or $m=2\ell$ and $n=\ell$ for some $\ell\in\Bbb Z$
(see Proposition \ref{1ref}).
As another interesting example,
the solution space of a system of linear equations
(familiar in elementary linear algebra)
is a reflection space in a vector space but not symmetric (see Example \ref{linearsol}).

The structure of symmetric spaces are simple.
Namely, a symmetric space of $G$ is just a union of some cosets,
$\bigcup_i (2S+s_i)$ in $S/2S$, for a subgroup $S$ of $G$.
However, reflection spaces are more complicated.
The author does not know the classification even for $G=\Bbb Z^n$
(or even for $G=\Bbb Z^2$).
We hope for someone to classify reflection spaces
(see Proposition \ref{refingeneral} and Example \ref{exref}).

We also study the reflection space generated by two elements
(see Definition \ref{refgensub}),
and prove some basic properties.
As an application,
we discuss some subalgebras of a multi-loop algebra
(see Example \ref{twogenerators} and \ref{pmtwogenerators}).

Eventually, we reach a nicer and simpler definition of Lie tori
as simply an $\frak R$-graded Lie algebra satisfying certain properties
(see Definition \ref{defnewL}),
where $\frak R$ is a {\bf locally extended affine root system}
defined in \cite{MY1}.
This new {\bf Lie $\frak R$-torus} can be identified with both a general Lie torus
and a normal Lie torus.
We also show that if $G$ is a torsion-free abelian group,
then any Lie $G$-torus is a Lie $\frak R$-torus.
Thus there is essentially no difference for them
 when $G$ is a torsion-free abelian group.
 One of the major benefits is that
 no abelian group $G$ is involved in the definition of a Lie $\frak R$-torus.
 So the description of a Lie $\frak R$-torus is much shorter,
 and we see that the definition only depends on 
 the locally extended affine root system $\frak R$.
 It is not difficult to see that
 if $\frak R$ is a finite irreducible root system,
 then
 the Lie $\frak R$-torus is nothing but a finite-dimensional split simple Lie algebra.
 If $\frak R$ is a locally finite irreducible root system (see \cite{LN1}),
 then
 the Lie $\frak R$-torus is a locally finite split simple Lie algebra,
 studied in \cite{NS} and \cite{St}.
 
 Let us explain more how nice Lie $\frak R$-tori are.
 If $\frak R$ is an affine root system
 defined in \cite{M},
 then
 the Lie $\frak R$-torus is a derived affine Lie algebra
 or a (twisted) loop algebra.
If $\frak R$ is an extended affine root system defined in \cite{S},
 then
 the Lie $\frak R$-torus is a central extension of the centerless core of an extended affine Lie algebra,
 studied in \cite{BGK}, \cite{BGKN} and \cite{AABGP}.
 Actually,  Allison and Gao started to research the centerless core
 as a double graded Lie algebra, which is the origin of  a Lie torus (see \cite{AG}).
 Since then,
 many people studied Lie tori, for example, in
 \cite{BY}, \cite{AY}, \cite{AFY1}, \cite{AFY2}, \cite{AB}, \cite{F}, \cite{Y3}, etc.
 
Moreover,
 J. Morita and the author studied in \cite{MY1}
 a generalization of locally finite split simple Lie algebras
 and extended affine Lie algebras.
 Thus, if $\frak R$ is a locally affine root system defined in \cite{Y4},
 then
 the Lie $\frak R$-torus is the core of a {\bf locally affine Lie algebra}
 or a {\bf locally (twisted) loop algebra} (see \cite{N} and \cite{MY2}).
If $\frak R$ is a locally extended affine root system,
 then
 the Lie $\frak R$-torus is a central extension of the centerless core of a locally extended affine Lie algebra
(see  \cite{MY1}).

\medskip

 The author thanks Professor Jun Morita for
thoughtful discussions and suggestions.
Also, the author would like to thank the referee for 
helpful comments and suggestions.

 \enskip

\enskip

\section{Basic concepts}
\label{}

Let $G= (G,+,0)$ be an arbitrary abelian group.
Let $\Delta$ be a {\bf locally finite irreducible root system} (see \cite{LN1}),
and we denote the Cartan integer by  
$$\langle\mu,\nu\rangle:=\frac{2(\mu,\nu)}{(\nu,\nu)}$$ for $\mu,\nu\in\Delta$,
and also let $\langle 0,\mu\rangle:=0$ for all $\mu\in\Delta$.
Recall that $\Delta$ is called {\bf reduced}
if $2\alpha\notin\Delta$
for all $\alpha\in\Delta$.
We define the subset
$$\Delta^{\red}:=
\{\alpha\in\Delta\mid\ \frac{1}{2}\alpha\notin\Delta\}
$$
 of $\Delta$,
 which is a reduced locally finite irreducible root system.
 Note that $\Delta=\Delta^{\red}$ if $\Delta$ is reduced.
We review the notion of a Lie $G$-torus
introduced in \cite{Ne},
which we call here
a normal locally Lie $G$-torus.
(Originally, it is defined for a finite irreducible root system $\Delta$,
but it is easily generalized to a locally finite irreducible root system.)

\begin{definition}\label{defgeneral} 
Let $\Delta$ be a locally finite irreducible root system.
A Lie algebra ${\mathcal L}$  is called a 
{\bf locally Lie 
$G$-torus of type $\Delta$}  if 
\begin{itemize}
\item[(LT1)]
${\mathcal L}$ has a decomposition into subspaces  
$${\mathcal L} =  \bigoplus_{\mu \in \Delta \cup \{0\},\  { g \in G}} \  {\mathcal L}^g_\mu$$ 
such that  $[{\mathcal L}^g_\mu, {\mathcal L}^h_\nu]  \subset {\mathcal L}^{g+h}_{\mu+\nu}$
for $\mu,\nu, \mu+\nu \in \Delta \cup \{0\}$
and $g,h\in G$;
\item[(LT2)]
 For every $g \in G$, \ \ ${\mathcal L}_0^g =  \sum_{\mu \in \Delta,\
h \in G} \  [{\mathcal L}_\mu^h,\  {\mathcal L}_{-\mu}^{g-h}]$; 
\item[(LT3)]
For each $0\neq x \in {\mathcal L}_\mu^g$ ($\mu \in \Delta,\  g \in G)$,  there exists 
$y \in {\mathcal L}_{-\mu}^{-g}$  so that 
$\mu^\vee : = [x, y] \in {\mathcal L}_0^0$ satisfies 
$[\mu^\vee , z] = \langle \nu,\mu\rangle z$ for all  $z \in {\mathcal L}_\nu^h$\
\ $(\nu \in \Delta \cup \{0\},\ h \in G)$;
\item[(LT4)]
$\langle\supp_G{\mathcal L}\rangle=G$,
where
$$\supp_G{\mathcal L} 
:= \{g \in G \mid {\mathcal L}_\mu^g \neq 0 \ \text{for some} \ \mu \in \Delta
\cup \{0\}\};$$
\item[(LT5)]
$\dim {\mathcal L}_\mu^g \leq 1$ for all $\mu \in \Delta$ and $g\in G$;
\item[(LT6)]
$\dim {\mathcal L}_\mu^0 = 1$ 
for all $\mu \in \Delta^\red$.
\end{itemize}

\remark
(i)  
Condition (LT4) is simply a convenience. If it fails to hold, we may replace
$G$ by the subgroup generated by $ \supp_G\mathcal L $. 

(ii)  It follows from (LT1) and (LT3) that 
$\mathcal L $ admits a grading by the root lattice $Q(\Delta)$:  if 
$$\mathcal L _\lambda := \bigoplus_{g \in G}\  \mathcal L _\lambda^g $$ for 
$\lambda \in Q(\Delta)$, where $\mathcal L _\lambda^g = 0$ if $\lambda \not \in \Delta \cup \{0\}$, then $\mathcal L  = \bigoplus_{\lambda \in Q(\Delta)}\ \mathcal L _\lambda$
and  $[\mathcal L _\lambda, \mathcal L _\mu] \subset \mathcal L _{\lambda+\mu}$. 

(iii)
$\mathcal L $ is also graded by the group $G$.  Namely, if
$$\mathcal L ^g := \bigoplus_{\mu \in \Delta \cup \{0\}}\ \mathcal L _\mu^g,$$ then
$\mathcal L  = \bigoplus_{g \in G}\ \mathcal L ^g$ and 
$[\mathcal L ^g,\mathcal L ^h] \subset \mathcal L ^{g+h}$.  

 \endremark

Now, we define a {\bf general locally Lie $G$-torus}
as a Lie algebra satisfying (LT1-5) above and instead of (LT6),
$$\text{(LT6)}'\ \ \ \ \
\mathcal L_\mu\neq 0\
\text{for all $\mu \in \Delta$}.$$
We often say a Lie torus 
when `general', `locally' or $G$ is clear from the context or no need to specify.
If $G\cong\Bbb Z^n$, then
${\mathcal L}$ is called a {\bf locally Lie $n$-torus} or simply a {\bf Lie $n$-torus}.
We call the rank of $\Delta$ the {\bf rank} of ${\mathcal L}$
and the type of $\Delta$ the {\bf type} of $\mathcal L$.
If $\mathcal L$ has the trivial center, then $\mathcal L$
is called {\bf centerless}.

Finally,
we say
a {\bf normal locally Lie $G$-torus}
when $\mathcal L$ satisfies  (LT6)$'$ and  (LT6),
not just (LT6).
By this, a Lie torus really has type $\text{BC}$
if $\Delta$ has type $\text{BC}$.
($\mathcal L$ might have reduced type even though $\Delta$ is nonreduced
if we only assume (LT6).)
We may omit the term `normal' or `general'
unless we compare two concepts.

\end{definition}

\medskip

Before giving examples of general Lie tori,
 we introduce basic concepts about 
 isomorphisms of graded algebras
 following [AF].
 \definition
 Let $A=\oplus_{g\in G}\ A_g$ 
and $A'=\oplus_{g'\in G'}\ A_{g'}$ 
be graded algebras,
where $G$ and $G'$ are groups.
An algebra isomorphism $\varphi:A\longrightarrow A'$
is called an {\bf isograded isomorphism} if
there exists a group isomorphism $\psi:G\longrightarrow G'$
such that $\varphi(A_g)\subset A'_{\psi(g)}$.
In particular, 
if $G=G'$ and $\psi$ is the identity map,
it is called a {\bf graded isomorphism}.
Also, we say that $A$ and $A'$ are
{\bf isograded isomorphic}
(resp. {\bf graded isomorphic})
if there exists an isograded isomorphism 
(resp. a graded isomorphism) between them.
We sometimes identify two $G$-graded algebras
if they are graded isomorphic.

A Lie $G$-torus is graded by 
the abelian group $\langle\Delta\rangle\times G$.
Thus
an isograded isomorphism between two Lie tori
means that they
are isograded isomorphic relative to such direct product groups.

For a group homomorphism
$s\in\hom(\langle\Delta\rangle, G)$ and a Lie $G$-torus $L$,
one can change the $G$-grading as follows, and
this new $G$-graded Lie algebra is called
an {\bf isotope} of $L$ by $s$, denoted by $L^{(s)}$:
$$(L^{(s)})_\alpha^{g}:=L_\alpha^{g+s(\alpha)}$$
for all $\alpha\in\langle\Delta\rangle$ and $g\in G$.

 \enddefinition

\lemma\label{isoreg}
An isotope $L^{(s)}$ of $L$ is in fact isograded isomorphic to $L$.
\endlemma
\proof
This follows from a general property for a direct product group.
Let $M$ and $G$ are abelian groups,
and let $s:M\longrightarrow G$ 
be a group homorphism.
Then
the map $f:M\times G\longrightarrow M\times G$
defined by $f(m,g)=(m,g+s(m))$
for $m\in M$ and $g\in G$
is a group automorphism of the product group $M\times G$.
In fact, $f$ is clearly a monomorphism.
For any $(m,g)\in M\times G$,
let $x=(m,g-s(m))$. Then $f(x)=(m,g)$, and so $f$ is onto.
\qed

\enskip

We will show that any isotope of a general Lie $G$-torus
is again a general Lie $G'$-torus,
where $G'$ is the subgroup of $G$
generated by $\supp_GL^{(s)}$.

 \definition
Let $\varphi$ be an isograded-isomorphism
from a Lie $G$-torus $L$ onto a Lie $G'$-torus $L'$,
and $\varphi$ is called a {\bf bi-isomorphism}
if the corresponding group isomorphism $f$ from
$\langle\Delta\rangle\times G$ onto
$\langle\Delta'\rangle\times G'$ decomposes
into group isomorphisms of each factor.
More precisely,
there exist 
a group isomorphism $w$ from
$\langle\Delta\rangle$ onto
$\langle\Delta'\rangle$ and
a group isomorphism $\psi$ from
$G$ onto
$G'$ such that $\varphi=w\times \psi$.

If a Lie $G$-torus $L$ is bi-isomorphic 
to an isotope of a Lie $G'$-torus $L'$, we say that $L$ is isotopic to $L'$,
denoted by $L\sim L'$.

 \enddefinition

\remark
Suppose that $\varphi$ is an isograded isomorphism from a $G$-graded algebra $A$
onto a $G'$-graded algebra.
In particular,
if $\varphi$ is the identity map on a $G$-graded algebra $A$
with a group isomorphism $\psi$
from $G$ onto $G'$,
we may say that $A$ is {\bf re-graded by $G'$ through $\psi$}.
For example,
an isotope $L^{(s)}$ of a Lie $G$-torus $L$ is re-graded 
(through $f$ in the proof of Lemma \ref{isoreg}).
Note that an isograded isomorphism of
Lie $G$-tori
is used for
$\langle\Delta\rangle\times G$-grading,
not just $G$-grading.
Also, an isotope $L^{(s)}$ does not change the degree of the first factor $\langle\Delta\rangle$.
However, we note that
$\supp L^{(s)}\neq\supp L$
and
$\supp_G L^{(s)}\neq\supp_G L$ 
in general.
Thus an isotope $L^{(s)}$ is not necessarily a Lie $G$-torus
since $P:=\supp_G L^{(s)}$ might be a proper subgroup of $G$.
But certainly $L^{(s)}$ is a Lie $P$-torus of the same type.
For such a case, 
{\bf we still say that $L$ is isotopic to a Lie $P$-torus $L^{(s)}$}.

\endremark

We will show that any general Lie torus 
can be re-graded to a normal Lie torus.

\enskip

\enskip

\section{Examples}
\label{}
 
Let us give some examples of general Lie tori.

\example\label{3exofnaoi}
Let $\{ e, f,h\}$ be an $\sll_2$-triple
so that $[e,f]=h$, $[h,e]=2e$ and $[h,f]=-2f$,
having the root system $\{\pm\alpha\}$ relative to $Fh$,
i.e.,
$\alpha$ is the linear form of $Fh$
such that $\alpha(h)=2$.

\medskip

(1)
Let $P_1:=F[t^{\pm p}]$
for an integer $p>1$ and
$L:=(e\otimes t^rP_1)\oplus (f\otimes t^{-r}P_1)
\oplus(h\otimes P_1)$,
where $r$ is a positive integer
so that $r$ and $p$ are coprime.
One can say that $L$ is the subalgebra of
a normal Lie $1$-torus 
$\sll_2(F[t^{\pm 1}])=\sll_2(F)\otimes F[t^{\pm 1}]$ 
generated by
$e\otimes t^r$, $f\otimes t^{-r}$
and $h\otimes t^{\pm p}$.
Then $L$ becomes a general Lie $\Bbb Z$-torus
by defining
$L_\alpha^{pm+r}=Fe\otimes t^{pm+r}$,
$L_{-\alpha}^{pm-r}=Ff\otimes t^{pm-r}$,
and $L_0^{pm}=Fh\otimes t^{pm}$
for all $m\in\Bbb Z$,
and all the other homogeneous spaces
$L_{\alpha}^{k}$,
$L_{-\alpha}^{k}$
and
$L_0^{k}$
are all $0$.
Let
$$S_\alpha:=\supp_\Bbb Z L_\alpha
=\{k\in\Bbb Z\mid L_\alpha^k\neq 0\}
\subset\supp_\Bbb Z L.$$ 
We see that
$S_{\alpha}=p\Bbb Z+r$,\ \  \ 
$S_{-\alpha}=p\Bbb Z-r$\ \ \ 
(so $S_{-\alpha}=-S_{\alpha}\neq S_{\alpha}$),\ \ \
and
$$\supp_\Bbb Z L=(p\Bbb Z+r)\cup(p\Bbb Z-r)\cup p\Bbb Z,$$
and hence $\langle\supp_\Bbb Z L\rangle=\Bbb Z$.

Let $s\in\hom(\langle\alpha\rangle,\Bbb Z)$
define by
$s(\alpha)=r$.
Then the $s$-isotope $L^{(s)}$ defined by
$$
(L^{(s)})_{\pm\alpha}^n:=L_{\pm\alpha}^{n\pm s(\alpha)}
=L_{\pm\alpha}^{n\pm r}
\quad\text{and}\quad
(L^{(s)})_{0}^n:=L_{0}^{n\pm s(0)}=L_{0}^{n}
$$
for all $n\in\Bbb Z$ can be identified with a normal Lie torus
since $(L^{(s)})_{\pm\alpha}^0\neq 0$.
In fact, first note that
$\supp_\Bbb Z (L^{(s)})=p\Bbb Z$, and then it is clear that
$L^{(s)}$ is a centerless normal Lie $p\Bbb Z$-torus.
Moreover,
one can easily check that
$L^{(s)}$ is graded isomorphic to $\sll_2(P_1)$.
Thus one can say that $L$ is isotopic to $\sll_2(P_1)$.

\begin{comment}
There is another way to check that $L$ is isograded isomorphic to
$\sll_2(P_1)$.
In fact, 
let 
$g=
\begin{pmatrix}
t^r&0\\
0&1
\end{pmatrix}
$.
Then $\varphi(x)=g^{-1}xg$ for $x\in L$ gives an isograded isomorphism
from $L$ onto $\sll_2(P_1)$.
Also, one can check that the group isomorphism determined by $\varphi$
is the same as the group isomorphism 
determined by the $s$-isotope $L^{(s)}$ above.
More precisely,
such a group isomorphism is given as
$f:\langle (\alpha, r), (0,p)\rangle\longrightarrow
\langle (\alpha, 0), (0,p)\rangle$
determined by
$f (\alpha, r)=(\alpha, 0)$
and $f (0, p)=(0, p)$.
\end{comment}

\medskip

(2)
Let $P_2:=F[t_1^{\pm p_1},t_2^{\pm p_2}]$ 
for some integers $p_1,p_2>1$,
and
$$L=
(e\otimes t_1^{r_1}t_2^{r_2}P_2)\oplus 
(f\otimes t_1^{-r_1}t_2^{-r_2}P_2)
\oplus(h\otimes P_2)$$
for some integers $r_1,r_2>1$
so that $(p_i, r_i)=1$ ($i=1,2$).
One can say that $L$ is
the subalgebra of
$\sll_2(F[t_1^{\pm 1},t_2^{\pm 1}])$ generated by
$e\otimes t_1^{r_1}t_2^{r_2}$, $f\otimes t_1^{-r_1}t_2^{-r_2}$,
$h\otimes t_1^{\pm p_1}$
and $h\otimes t_2^{\pm p_2}$.
Then $L$ is a general Lie $\Bbb Z^2$-torus,
and consider
the $s$-isotope $L^{(s)}$, where 
$s\in\hom(\langle\alpha\rangle,\Bbb Z^2)$
is defined as
$s(\alpha)=(r_1,r_2)$.
Note that
$$\supp_{\Bbb Z\times\Bbb Z} L_\alpha
=
(p_1\Bbb Z+r_1)\times (p_2\Bbb Z+r_2)
\quad
\text{and}
\quad
\supp_{\Bbb Z\times\Bbb Z}  L^{(s)}_\alpha
=p_1\Bbb Z\times p_2\Bbb Z,
$$
and one can show that
$L^{(s)}$
is a normal Lie $2$-torus.
Thus
$L$ is isotopic to 
a normal Lie $2$-torus.

\medskip

(3)
Let $\{ e_i, f_i\mid i=1,2\}$ be a set of
Chevalley generators of $\sll_3(F)$
with a Cartan subalgebra $\mathfrak h=\spa \{[e_i, f_i]\mid i=1,2\}$
and the root system
$\{\pm\alpha_1,\pm\alpha_2,\pm(\alpha_1+\alpha_2)\}$.
Let $L$ be the subalgebra of 
$\sll_3(F[t^{\pm 1}])$ generated by
$e_i\otimes t^{r_i}$, $f_i\otimes t^{-r_i}$
and $\frak h\otimes t^{\pm p}$.
Then $L$ is a general Lie $1$-torus.
Note that $S_{\alpha_i}=p\Bbb Z+r_i$
and $S_{\alpha_1+\alpha_2}=p\Bbb Z+r_1+r_2$
can be all different sets.

Consider the $s$-isotope $L^{(s)}$, where 
$s\in\hom(\langle\alpha_1,\alpha_2\rangle,\Bbb Z)$
is defined as
$s(\alpha_i)=r_i$.
Note that
$\supp_{\Bbb Z}  L=\Bbb Z
\supset p\Bbb Z
=\supp_{\Bbb Z}  \big(\sll_3(P_1)\big)$,
and one can show that
$L^{(s)}$
is  a normal Lie $p\Bbb Z$-torus
which is graded isomorphic to $\sll_3(P_1)$.
Thus $L$ is isotopic to $\sll_3(P_1)$.

\medskip

(4)
Next examples are
twisted loop algebras 
which look different.
Let
$\frak I$ be an arbitrary index set (possibly infinite)
with $|\frak I|\geq 2$.
Let $V$ be a $|\frak I|$-dimensional vector space over $\Bbb Q$
with a positive definite symmetric bilinear form.
Let
$\{\epsilon_i\mid i\in\frak I\}$
be an orthonormal basis of $V$.
Let 
$$\text D_\frak I=\{\pm\epsilon_i\pm\epsilon_j\mid i,j\in\frak I,\ i\neq j\}
\subset
\text C_\frak I
=\{\pm\epsilon_i\pm\epsilon_j,\pm 2\epsilon_i\mid i,j\in\frak I,\ i\neq j\}$$
be locally finite irreducible root systems
of type $\text D_\frak I$ and $\text C_\frak I$
(see \cite{LN1} or \cite{NS}).
Let
\begin{align*}
\text D_\frak I'&:=\{\pm(\epsilon_i-\epsilon_j)\mid i,j\in\frak I,\ i\neq j\}\\
\text D_\frak I''&:=\{\pm(\epsilon_i+\epsilon_j)\mid i,j\in\frak I,\ i\neq j\}\\
\text C_\frak I'&:=\{\pm(\epsilon_i+\epsilon_j)\mid i,j\in\frak I\}
\quad\text{so that}\\
\text D_\frak I&=\text D_\frak I'\cup\text D_\frak I''
\quad\text{and}\quad
\text C_\frak I=\text D_\frak I'\cup\text C_\frak I'.
\end{align*}
Let 
$$s_1=
\begin{pmatrix}
0&\iota\\
\iota&0
\end{pmatrix}
\quad\text{and}\quad
s_2=
\begin{pmatrix}
0&-\iota\\
\iota&0
\end{pmatrix}
$$
be the matrices of size $2\frak I$,
where
$\iota$ is the identity matrix of size $\frak I$.
Define the automorphisms $\sigma_i$ ($i=1, 2$) of
$\sll_{2\frak I}(F)$ by
$$
\sigma_i(x)=s_i^{-1}x^Ts_i
$$
for $x\in\sll_{2\frak I}(F)$,
where $x^T$ is the transpose of $x$.
Then the fixed algebra $\sll_{2\frak I}(F)^{\sigma_i}$
of $\sll_{2\frak I}(F)$
is a
locally finite split simple Lie algebra of type 
$\text D_\frak I$ or $\text C_\frak I$.
Thus,
let 
$$\sll_{2\frak I}(F)^{\sigma_1}=\frak g^D=\frak h\oplus\bigoplus_{\xi\in\text D_\frak I}\frak g_\xi^D
\quad\text{and}\quad
\sll_{2\frak I}(F)^{\sigma_2}=\frak g^C=\frak h\oplus\bigoplus_{\xi\in\text C_\frak I}\frak g_\xi^C.$$
Note that one can take the same Cartan subalgebra $\frak h$
for the types $\text D_\frak I$ and $\text C_\frak I$,
and also $\frak g_\xi^D=\frak g_\xi^C$
for all $\xi\in D'_\frak I$.
We extend $\sigma_i$ to the loop algebra
$\sll_{2\frak I}(F)\otimes F[t^{\pm 1}]$,
denoted $\hat\sigma_i$,
as
$$
\hat\sigma_i(x\otimes t^m)=(-1)^m\sigma_i(x)\otimes t^m
$$
for $x\in\sll_{2\frak I}(F)$ and $m\in\Bbb Z$.
Let
$$
T(D):=
(\sll_{2\frak I}(F)\otimes F[t^{\pm 1}])^{\hat\sigma_1}
\quad\text{and}\quad
T(C):=
(\sll_{2\frak I}(F)\otimes F[t^{\pm 1}])^{\hat\sigma_2}
$$
be the fixed algebras, which are usually called
the twisted loop algebras by $\hat\sigma_i$.
Let
$$V_\xi^i=\{x\in\sll_{2\frak I}(F)\mid
[h, x]=\xi(h)x\ \text{for all $h\in\frak h$ and $\sigma_i(x)=-x$}\}.$$
Then we can show that
\begin{align*}
V_\xi^1&=V_\xi^2\quad\text{for all $\xi\in\text D_\frak I'\cup\{0\}$}\\
V_\xi^1&=\frak g_\xi^C\quad\text{for all $\xi\in\text C_\frak I'$}\\
V_\xi^2&=\frak g_\xi^D\quad\text{for all $\xi\in\text D_\frak I''$},
\end{align*}
and moreover, letting
$V_\xi:=V_\xi^1$ for all $\xi\in\text D_\frak I'\cup\{0\}$,
we obtain
$$
T(D)
=(\frak h\oplus\bigoplus_{\xi\in\text D_\frak I}\frak g_\xi^D)
\otimes F[t^{\pm 2}]
\oplus
(V_0\oplus\bigoplus_{\xi\in\text D_\frak I'}V_\xi
\oplus\bigoplus_{\xi\in\text C_\frak I'}\frak g_\xi^C)
\otimes tF[t^{\pm 2}]
$$
and
$$
T(C)
=(\frak h\oplus\bigoplus_{\xi\in\text C_\frak I}\frak g_\xi^C)
\otimes F[t^{\pm 2}]
\oplus
(V_0\oplus\bigoplus_{\xi\in\text D_\frak I'}V_\xi
\oplus\bigoplus_{\xi\in\text D_\frak I''}\frak g_\xi^D)
\otimes tF[t^{\pm 2}].
$$
The latter algebra $T(C)$ is a so-called twisted loop algebra of type $\text C_\ell^{(2)}$
(or $\text A_{2\ell+1}^{(2)}$ in Kac label)
when $|\frak I|=\ell$ is finite.
What is the former algebra $T(D)$ then?

It seems that $T(D)$ dose not appear on the list of Kac-Moody Lie algebras
(see e.g. \cite{K}).
We can at leaast check that $T(D)$ 
is a general locally Lie $1$-torus of type $\text C_\frak I$
(but not 
$\text D_\frak I$,
see the axiom (LT1) of a Lie $G$-torus).
Of course,
$T(C)$
is a normal locally Lie $1$-torus of type $\text C_\frak I$.

\medskip

We can now answer the question (see also \cite{H}).

\begin{proposition}
$T(D)$ is an isotope of $T(C)$.
\end{proposition}
\proof
Define the group homomorphism 
$s:\langle\text C_\frak I\rangle\longrightarrow\Bbb Z$
by
$s(\epsilon_{i_0}-\epsilon_i)=0$ and
$s(2\epsilon_{i_0})=1$
for a fixed $i_0\in\frak I$ and all $i\in\frak I$.
Then 
$$
s(\epsilon_i+\epsilon_{i_0})
=s(\epsilon_i-\epsilon_{i_0}+2\epsilon_{i_0})
=s(\epsilon_i-\epsilon_{i_0})+s(2\epsilon_{i_0})=0+1=1,
\quad\text{and hence}$$
$$
s(2\epsilon_i)
=s(\epsilon_i-\epsilon_{i_0}+\epsilon_i+\epsilon_{i_0})
=s(\epsilon_i-\epsilon_{i_0})+s(\epsilon_i+\epsilon_{i_0})=0+1=1.$$
Also, we have
$$s(\epsilon_i-\epsilon_j)
=s(\epsilon_i-\epsilon_{i_0}+\epsilon_{i_0}-\epsilon_j)
=s(\epsilon_i-\epsilon_{i_0})+s(\epsilon_{i_0}-\epsilon_j)=0+0=0
\quad\text{and}$$
$$s(\epsilon_i+\epsilon_j)
=s(\epsilon_i+\epsilon_{i_0}+\epsilon_j-\epsilon_{i_0})
=s(\epsilon_i+\epsilon_{i_0})+s(\epsilon_j-\epsilon_{i_0})=1+0=1.$$
Let 
$$T(D)=\bigoplus_{(\alpha,m)\in
(\text C_\frak I\cup\{0\})\times\Bbb Z} T_\alpha^m,
\quad\text{and}\quad
P_\alpha^m:=T_\alpha^{m+s(\alpha)}.$$
Then we have
\begin{align*}
P_{\epsilon_i-\epsilon_j}^{2m}&
=T_{\epsilon_i-\epsilon_j}^{2m+s(\epsilon_i-\epsilon_j)}
=T_{\epsilon_i-\epsilon_j}^{2m}
=\frak g_{\epsilon_i-\epsilon_j}^D\otimes t^{2m}\\
P_{\epsilon_i-\epsilon_j}^{2m-1}
&=T_{\epsilon_i-\epsilon_j}^{2m-1+s(\epsilon_i-\epsilon_j)}
=T_{\epsilon_i-\epsilon_j}^{2m-1}
=V_{\epsilon_i-\epsilon_j}\otimes t^{2m-1}\\
P_{\epsilon_i+\epsilon_j}^{2m}&
=T_{\epsilon_i+\epsilon_j}^{2m+s(\epsilon_i+\epsilon_j)}
=T_{\epsilon_i+\epsilon_j}^{2m+1}
=\frak g_{\epsilon_i+\epsilon_j}^C\otimes t^{2m+1}
=(\frak g_{\epsilon_i+\epsilon_j}^C\otimes t)\otimes t^{2m}\\
P_{\epsilon_i+\epsilon_j}^{2m-1}
&=T_{\epsilon_i+\epsilon_j}^{2m-1+s(\epsilon_i+\epsilon_j)}
=T_{\epsilon_i+\epsilon_j}^{2m}
=\frak g_{\epsilon_i+\epsilon_j}^D\otimes t^{2m}
=(\frak g_{\epsilon_i+\epsilon_j}^D\otimes t)\otimes t^{2m-1}\\
P_{2\epsilon_i}^{2m}&
=T_{2\epsilon_i}^{2m+s(2\epsilon_i)}
=T_{2\epsilon_i}^{2m+1}
=\frak g_{2\epsilon_i}^C\otimes t^{2m+1}
=(\frak g_{2\epsilon_i}^C\otimes t)\otimes t^{2m}\\
P_{2\epsilon_i}^{2m-1}&
=T_{2\epsilon_i}^{2m-1+s(2\epsilon_i)}
=T_{2\epsilon_i}^{2m}
=0\\
P_{0}^{2m}&
=T_{0}^{2m+s(0)}
=T_{0}^{2m}
=\frak h\otimes t^{2m}\\
P_{0}^{2m-1}&
=T_{0}^{2m-1+s(0)}
=T_{0}^{2m-1}
=V_0\otimes t^{2m-1},\\
\end{align*}
and in particular, the subalgebra $P^0$ of $T(D)$
has the following decomposition:
\begin{align*}
P^0=
\bigoplus_{\xi\in\text C_\frak I\cup\{0\}}P_\xi^0
&=\bigoplus_{\xi\in\text D_\frak I'}\frak g_\xi^D\oplus
\bigoplus_{\xi\in\text C_\frak I'^+}
\frak g_{\pm\xi}^C\otimes t^{\pm 1}\oplus
\frak h,
\end{align*}
where
$\text C_\frak I'^+=\{\epsilon_i+\epsilon_j\mid i,j\in\frak I\}$.
We can see that $P^0$
is isomorphic to $\frak g^C$
through 
$$\frak g_{\pm\xi}^C\otimes t^{\pm 1}\ni x_{\pm}\otimes t^{\pm 1}\ \mapsto\ x_{\pm}\in\frak g_{\pm\xi}^C$$
for $\xi\in\text C_\frak I'^+$
and the identity map for the other root spaces.
Moreover, 
\begin{equation*}\label{tdgrading}
T(D)^{(s)}=\bigoplus_{(\alpha,m)\in
(\text C_\frak I\cup\{0\})\times\Bbb Z} P_\alpha^m
\end{equation*}
is graded isomorphic to $T(C)$
through 
\begin{align*}
P_{\pm(\epsilon_i+\epsilon_j)}^{2m}=\frak g_{\pm(\epsilon_i+\epsilon_j)}^C\otimes t^{2m\pm1}
\ni x_{\pm}\otimes t^{2m\pm1}
\ &\mapsto\ 
x_{\pm}\otimes t^{2m}
%\in\frak g_{\pm(\epsilon_i+\epsilon_j)}^C\otimes t^{2m}
\in T(C)_{\pm(\epsilon_i+\epsilon_j)}^{2m}
\\
P_{\pm(\epsilon_i+\epsilon_j)}^{2m\mp 1}=\frak g_{\pm(\epsilon_i+\epsilon_j)}^D\otimes t^{2m}
\ni x_{\pm}\otimes t^{2m}
\ &\mapsto\ 
x_{\pm}\otimes t^{2m\mp 1}
%\in\frak g_{\pm(\epsilon_i+\epsilon_j)}^D\otimes t^{2m\mp 1}
\in T(C)_{\pm(\epsilon_i+\epsilon_j)}^{2m\mp 1}
\\
P_{\pm 2\epsilon_i}^{2m}=\frak g_{\pm 2\epsilon_i}^C\otimes t^{2m\pm 1}
\ni x_{\pm}\otimes t^{2m\pm 1}
\ &\mapsto\ 
x_{\pm}\otimes t^{2m}
%\in\frak g_{\pm 2\epsilon_i}^C\otimes t^{2m}
\in T(C)_{\pm 2\epsilon_i}^{2m}
\end{align*}
and the identity map for the other homogeneous spaces.
Thus 
$T(D)$ is isotopic to $T(C)$.
(In particular,
$T(D)^{(s)}$ is
a normal locally Lie torus of type $\text C_\frak I$.)
%with grading pair $(P^0,\frak h)$.
\qed

\endexample

\enskip

\section{Relation between general and normal Lie tori}
\label{}

In Example \ref{3exofnaoi} 
we learned how to get a normal Lie torus
from a general Lie torus.
The process can be generalized.
Let us first recall reflectable sections (see \cite{Y4}).

Let $(\Delta,V)$ be a locally finite irreducible root system.
We define $\sigma_\alpha$ for $\alpha\in\Delta$ by
$$
\sigma_\alpha(\beta)=\beta-\langle \beta,\alpha\rangle\alpha
$$
for $\beta\in\Delta$.

A basis $\Pi$ of $V$ as a vector space
is called a {\bf reflectable base} of $\Delta$
if $\Pi\subset\Delta$ and for any $\alpha\in\Delta^\red$,
$$\alpha=\sigma_{\alpha_1}\cdots\sigma_{\alpha_k}(\alpha_{k+1})$$
for some $\alpha_1, \ldots, \alpha_{k+1}\in\Pi$.
(Any root can be obtained by reflecting a root of $\Pi$ relative to the hyperplanes determined by $\Pi$.)
This is a well-known property of a root base in a reduced finite root system.
It is known that a locally finite irreducible root system which is countable has a root base,
but this is not the case for uncountable ones
(see [LN1, \S 6]).
However, it is proved in [LN2, Lem.5.1]
that there exists 
a reflectable base 
in a reduced locally finite irreducible root system even if it is uncountable.

\definition\label{grauto}
Let $\Bbb Z^{(\frak I)}$ 
be a free abelian group of rank $\frak I$
and let
$G$ 
be an arbitrary abelian group.
Let
$\frak B=\{\mu_i\}_{i\in\frak I}$
be
a basis of $\Bbb Z^{(\frak I)}$.
Fix some
$g_i\in G$
for $i\in\frak I$.
Then the group homomorphism $s\in\hom(\Bbb Z^{(\frak I)}, G)$
defined by
$s(\mu_i)=g_i$
for all $i$
is called the {\bf shift}
relative to $\frak B$ and $\{g_i\}_{i\in\frak I}$.

\enddefinition

\begin{comment}
\lemma\label{nouse}
 Let $\Delta$ be a locally finite root system,
and let
 $s\in\hom(\langle\Delta\rangle, G)$.
 For $\alpha\in\Delta$,
 suppose 
 $\alpha=\sigma_{\alpha_1}\cdots\sigma_{\alpha_k}(\alpha_{k+1})$
for some $\alpha_1, \ldots, \alpha_{k+1}\in\Delta$.
Then
 $$s(\alpha)=\tilde\sigma_{\alpha_1}\cdots\tilde\sigma_{\alpha_k}\big(s(\alpha_{k+1})\big),$$
where
 $$\tilde\sigma_{\alpha_i}\big(s(\alpha_{j})\big)
 :=s(\alpha_{j})-\langle\alpha_j,\alpha_i\rangle s(\alpha_i).$$
 \endlemma
\proof
We use induction on $k$.
For $k=1$, i.e.,
if $\alpha=\sigma_{\alpha_1}(\alpha_2)$,
then
$$s\big(\sigma_{\alpha_1}(\alpha_2)\big)
=s\big(\alpha_2-\langle\alpha_2,\alpha_1\rangle\alpha_1\big)
=s(\alpha_{2})-\langle\alpha_2,\alpha_1\rangle s(\alpha_1)
=\tilde\sigma_{\alpha_1}\big(s(\alpha_{2})\big).
$$
Suppose that 
 $\beta=\sigma_{\alpha_2}\cdots\sigma_{\alpha_k}(\alpha_{k+1})$
for some $\alpha_2, \ldots, \alpha_{k+1}\in\Delta$
implies
 $$s(\beta)=\tilde\sigma_{\alpha_2}\cdots\tilde\sigma_{\alpha_k}\big(s(\alpha_{k+1})\big).$$
Then we obtain
$$s\big(\sigma_{\alpha_1}(\beta)\big)
=s(\beta)-\langle\beta,\alpha_1\rangle s(\alpha_1)
=\tilde\sigma_{\alpha_1}\big(s(\beta)\big)
=\tilde\sigma_{\alpha_1}\cdots\tilde\sigma_{\alpha_k}\big(s(\alpha_{k+1})\big).
\qed$$

(This lemma is not used in the paper,
but it may be useful if one tries to classify the root systems extended by $G$
for general Lie $G$-tori.)
\end{comment}

\lemma\label{usereffortori}
 Let $\Delta$ be a locally finite irreducible root system,
and let
$${\mathcal L} =  \bigoplus_{\mu \in \Delta \cup \{0\},\  { g \in G}} \  {\mathcal L}^g_\mu$$ 
be a general Lie $G$-torus.
Then, for  $s\in\hom(\langle\Delta\rangle, G)$,
we have
 $$
 L_{\alpha}^{s(\alpha)}\neq 0\ \
 \text{and}\ \
 L_{\beta}^{s(\beta)}\neq 0\ \
 \Longrightarrow\ \
 L_{\sigma_{\alpha}(\beta)}^{s(\sigma_{\alpha}(\beta))}
 \neq 0,
$$
and moreover,
 $$
 L_{\alpha_1}^{s(\alpha_1)}\neq 0,\  \ldots,\
L_{\alpha_k}^{s(\alpha_k)}\neq 0\ \
 \text{and}\ \
 L_{\alpha_{k+1}}^{s(\alpha_{k+1})}\neq 0\ \
 \Longrightarrow\ \
 L_{\alpha}^{s(\alpha)}
 \neq 0,
$$
where
 $\alpha=\sigma_{\alpha_1}\cdots\sigma_{\alpha_k}(\alpha_{k+1})$.
\endlemma
\proof
For $(\alpha, g), (\beta, g')\in\Delta\times G$,
we define
$$
\sigma_{(\alpha,g)}\big((\beta, g')\big)
:=\big(\sigma_\alpha(\beta),\ g'-\langle\beta,\alpha\rangle g\big).
$$
Suppose that there exists $0\neq x\in L_\alpha^g$.
Take $y\in L_{-\alpha}^{-g}$
such that $\{x, y, [x,y]\}$ is an $\sll_2$-triple,
and let
$\theta_\alpha^g:=\exp(\ad x)\exp(-\ad y)\exp(\ad x)$.
Then
$\theta_\alpha^g$ is an automorphism of $L$,
satisfying that
$\theta_\alpha^g(L_\beta^{g'})=L_{\sigma_\alpha(\beta)}^{g'-\langle\beta,\alpha\rangle g}$.
(One can prove this by the same way as in \cite[Prop.1.27]{AABGP}.)
In particular,
$L_\alpha^g\neq 0$ and $L_\beta^{g'}\neq 0$ implies that
$L_{\sigma_\alpha(\beta)}^{g'-\langle\beta,\alpha\rangle g}\neq 0$.
Thus
we have
 $$
 L_{\alpha}^{s(\alpha)}\neq 0\ \
 \text{and}\ \
 L_{\beta}^{s(\beta)}\neq 0\ \
 \Longrightarrow\ \
 L_{\sigma_{\alpha}(\beta)}^{s(\sigma_{\alpha}(\beta))}=
 L_{\sigma_{\alpha}(\beta)}^{s(\beta)-\langle\beta,\alpha\rangle s(\alpha)}
 \neq 0.
$$
So, for the last assertion,
it is true for $k=1$.
But if $L_{\beta}^{s(\beta)}
 \neq 0$,
where
 $\beta=\sigma_{\alpha_2}\cdots\sigma_{\alpha_k}(\alpha_{k+1})$,
then we have
$0\neq L_{\sigma_{\alpha_1}(\beta)}^{s(\beta)-\langle\beta,\alpha\rangle s(\alpha_1)}
=L_{\sigma_{\alpha_1}(\beta)}^{s(\sigma_{\alpha_1}(\beta))}=L_{\alpha}^{s(\alpha)}$.
\qed

\medskip

 \theorem\label{isotopicth}
Any general locally Lie $G$-torus $L$ of type $\Delta$
is isotopic to a normal Lie $P$-torus,
where $P$ is a subgroup of $G$.
  \endtheorem
  \proof    
  Let $\Pi=\{\mu_i\}_{i\in\frak I}$
 be a reflectable base of $\Delta$.
  Then there exist $g_i\in G$ for all $i\in\frak I$
  so that $L_{{\bf \mu}_i}^{g_i}\neq 0$ for all $i\in\frak I$
  by the axiom (LT6)$'$ of a general Lie $G$-torus.
Using the shift
$s\in\hom(\Bbb Z^{(\frak I)}, G)$
defined by
$s(\mu_i)=g_i$
for all $i\in\frak I$,
one gets the $s$-isotope $L^{(s)}$.
Now, we have
$\Pi\times 0\subset
\supp_{\langle\Delta\rangle\times G}L^{(s)}$
since $(L^{(s)})_{\mu_i}^0=L_{\mu_i}^{s(\mu_i)}=L_{{\bf \mu}_i}^{g_i}\neq 0$.
Then since $\Pi$ is reflectable,
we get 
$(L^{(s)})_{\alpha}^0=L_{\alpha}^{s(\alpha)}\neq 0$ for all $\alpha\in\Delta^\red$,
by Lemma \ref{usereffortori}.
Next, for $(\alpha,g)\in\langle\Delta\rangle\times G$,
we have
$\dim(L^{(s)})_{\alpha}^g
=\dim L_\alpha^{g+s(\alpha)}\leq 1$.
Hence we have the $1$-dimensionality (LT5).
Also, since $s(-\alpha)=-s(\alpha)$,
(LT3) holds.
Finally, 
since $L_0$ and $L_\alpha$ are unchanged by taking isotopes,
 the property
$(L^{(s)})_0=\sum_{\alpha\in\Delta}[(L^{(s)})_\alpha, (L^{(s)})_{-\alpha}]$
holds.
Thus let $P$ be the subgroup of $G$ generated by 
$\supp_G L^{(s)}$.
Then
$L^{(s)}$ is a normal $P$-torus,
and so
$L$ is isotopic to the normal $P$-torus.
    \qed

\medskip

\begin{remark}
A locally Lie $n$-torus $L$ is always a free module over the centroid
(see \cite{BN}).
Thus by Theorem \ref{isotopicth},
a general locally Lie $n$-torus has the same property.
We call the {\bf central rank} of $L$ the rank of 
the free module $L$ over the centroid.

\end{remark}

\definition
We call a symmetric invariant bilinear form on a Lie algebra $L$
simply a {\it form}. 
Here `invariant' is in the sense that $([x,y],z)=(x,[y,z])$
for all $x,y,z\in L$.
Note that if 
a $\Delta$-graded Lie algebra $L=\oplus_{\mu\in\Delta\cup\{0\}}\ L_{\mu}$
has a form $(\cdot,\cdot)$,
then 
$(L_\mu,L_\nu)=0$ unless $\mu+\nu= 0$ for $\mu,\nu\in\Delta\cup\{0\}$.

Moreover, if $L$ is $G$-graded and a form satisfies the property that
$(L^g,L^k)=0$ unless $g+k=0$ for all $g,k\in G$,
then the form is called a {\bf graded form}.
\enddefinition

The existence of a graded form on a Lie $n$-torus
is shown in \cite{Y3}.
Thus we have:

 \corollary\label{efrom}
 Any general Lie $n$-torus admits a nonzero graded form.
  \endcorollary
  \proof
  It follows from the implication
 $\mu+\nu=0$
  $\Longrightarrow$
 $s(\mu)+s(\nu)=0$.
  \qed

\medskip

We do not know the existence of a nonzero graded form
for a general locally Lie $G$-torus,
but if $G$ is torsion-free, 
it will be affirmative.

\enskip

\enskip

\section{Naoi tori}
\label{}

We introduce a Naoi torus defined in \cite{Na}.
(We slightly modified for our convenience.)
  
\definition\label{defNaoi}
A $\Bbb Z^n$-graded Lie algebra $L=\oplus_{g\in \Bbb Z^n}\ L^g$
is called a {\it Naoi torus} if the following conditions are satisfied:

\item[(N1)]
$L$ is graded simple.
\item[(N2)]
The central grading group has rank $n$.
\item[(N3)]
$L$ admits a nondegenerate graded form $(\ \ , \ \ )$.
\item[(N4)]
There exists an ad-diagonalizable subalgebra $\mathfrak h\subset L^0$
($\frak h$ is automatically abelian)
such that 
$$L^0\cap L_0=\frak h$$ 
and the set of roots
$$\Delta:=
\{0\neq\mu\in\frak h^*\mid L_\mu\neq 0\},
$$
where
$L_{\mu}=\{x\in L\mid [h,x]\in\mu(h)x\ \text{for all}\ h\in\mathfrak h\}$,
forms a finite irreducible root system
in the vector space $\spa_\Bbb Q \Delta$ spanned by $\Delta$ over $\Bbb Q$
relative to the induced form by scaling the above graded form on $L$.

More precisely,
by the property $L^0\cap L_0=\frak h$,
the restriction of the nondegenerate graded form on $\frak h$ is still nondegenerate.
Thus for $\mu\in\frak h^*$,
one can define
 a unique element  $t_\mu$ in $\frak h$
so that $(t_\mu,h)=\mu(h)$ for all $h\in\frak h$.
(Note that $t_0=0$.)
Then one can define a symmetric bilinear form
on $\spa_\Bbb Q \Delta$ as
$(\mu,\nu):=(t_\mu,t_\nu)$ for $\mu,\nu\in\frak h^*$.
Thus the latter part of (N4) says that after scaling the graded form,
the symmetric bilinear form on $\spa_\Bbb Q \Delta$ becomes positive definite
on $\Bbb Q$, and $\Delta$ is a finite irreducible root system in $\spa_\Bbb Q \Delta$
relative to this positive definite form.
In particular, $\frak h$ is finite-dimensional.
\enddefinition

Any centerless normal Lie $n$-torus of finite rank is clearly a Naoi torus.
Also, any of the Lie algebras in Example  \ref{3exofnaoi} is a Naoi torus.

\lemma\label{likeEALA}
A Naoi torus $L=\bigoplus_{g\in \Bbb Z^n}\ L^g$ 
has 
the docompostion
$L^g=\oplus_{\mu\in\Delta\cup\{0\}}\ L_{\mu}^g$,
where $L_{\mu}^g=L_{\mu}\cap L^g$.
In particular
$L$ has
the double grading
$$L=\bigoplus_{\mu\in\Delta\cup\{0\}}\ \bigoplus_{g\in \Bbb Z^n}\ L_{\mu}^g,
$$
and $L_{\mu}=\bigoplus_{g\in \Bbb Z^n}\ L_{\mu}^g$.
Moreover, we have
$$
[x,y]=(x,y)t_\mu
$$
for $x\in L_\mu^g$ and $y\in L_{-\mu}^{-g}$.
In particular,
$[x,y]=0$
for $x\in L_0^g$ and $y\in L_0^{-g}$.

\endlemma

\proof
The first assertion is clear since each $L^g$ is an $\frak h$-weight module.
For the next assertion,
we have  $([x,y]-(x,y)t_\mu,h)=([x,y],h)-(x,y)(t_\mu,h)
=(x,[y,h])-\mu(h)(x,y)=0$
for all $h\in\frak h$.
Hence $[x,y]=(x,y)t_\mu$.
The last assertion follows from $t_0=0$.
\qed

 \theorem\label{naoiisgen}
 Any Naoi torus is a centerless general Lie $n$-torus of finite central rank.
 
 Conversely, any centerless general Lie $n$-torus of finite central rank
 is a Naoi torus.
  \endtheorem
  \proof
  Let $L$ be a Naoi torus. 
  Then,
 by Lemma \ref{likeEALA},
$$L'=\bigoplus_{\mu\in\Delta}\ L_{\mu}\oplus\sum_{\mu\in\Delta}[L_\mu,L_{-\mu}]$$
is a graded ideal of $L$, and hence
we get $L=L'$.
In particular,
$L_0=\sum_{\mu\in\Delta}[L_\mu,L_{-\mu}]$.
Thus (LT2) holds.
Next,
for $0\neq x\in L_\mu^g$,
let 
$y$ be an element in $L_{-\mu}^{-g}$
such that $(x,y)=\frac{2}{(\mu,\mu)}$.
Then  by Lemma \ref{likeEALA},
letting
$\mu^\vee:=[x,y]=\frac{2t_\mu}{(\mu,\mu)}$
and for $z\in L_{\nu}^{k}$,
we have
$[\mu^\vee, z]=\nu(\mu^\vee)z=\frac{2\nu(t_\mu)}{(\mu,\mu)}z
=\frac{2(\nu,\mu)}{(\mu,\mu)}z
=\langle\nu,\mu\rangle z$,
and hence (LT3) holds.
For the 1-dimesionality (LT5),
we first claim that for $0\neq y\in L_{-\mu}^{-g}$, the linear map
$$
\ad y:L_\mu^g\longrightarrow Ft_\mu
$$
is injective.
In fact, suppose $(x',y)t_\mu=0$ for $x'\in L_\mu^g$.
Then we have $(x',y)=0$,
and so $[y,x']=0$.
But note that $\{x,y,\mu^\vee\}$ is an $\sll_2$-triple,
and $[\mu^\vee,x']=2x'$.
So the identity $[y,x']=0$ cannot happen by the $\sll_2$-theory
unless $x'=0$.
Hence our claim is settled.
Then we have $\dim L_\mu^g\leq\dim Ft_\mu=1$.
Finally, note that the rank of the central grading group of a Naoi torus
is equal to the rank of the grading group.
Hence the central rank of a Naoi torus is finite.
  
  The converse is clear since 
  a general Lie $n$-torus is isotopic to a normal Lie torus
  (see also Corollary \ref{efrom}).
  \qed

\enskip

\enskip

\section{Reflection spaces}
\label{}

Let
$$S_\alpha:=\supp_GL_\alpha 
=\{g\in G\mid L_\alpha^g\neq 0\}
\subset\supp_G L$$ 
for a locally general Lie $G$-torus $L$
of type $\Delta$.
This set for the normal case was classified in \cite{Y2},
but the set is quite wild compared with the normal case.

First since the reflections by the locally extended affine root system
via $(\alpha, s)\leftrightarrow \alpha+s$,
acts on $\supp L$,
we have, for $s\in S_\alpha$,
$$
\sigma_{\alpha+s}(\alpha+s)
=\alpha+s-\langle\alpha+s,\alpha+s\rangle(\alpha+s)
=-\alpha-s\in\supp L.
$$
Thus $-s\in S_{-\alpha}$,
and so $-S_{\alpha}\subset  S_{-\alpha}$.
Similarly, we have
$-S_{-\alpha}\subset  S_{\alpha}$, and hence
\begin{equation}\label{weakrel}
-S_\alpha= S_{-\alpha}.
\end{equation}
Also,
we have
$$
\sigma_{\alpha+t}(\alpha+s)
=\alpha+s-\langle\alpha+s,\alpha+t\rangle(\alpha+t)
=-\alpha+s-2t\in\supp L,
$$
and hence $s-2t\in S_{-\alpha}$
for all $s,t\in S_\alpha$.
Thus 
$S_\alpha-2S_\alpha\subset S_{-\alpha}$,
and
by \eqref{weakrel}, 
we get
\begin{equation}\label{weakrefs}
2S_\alpha-S_\alpha\subset S_{\alpha}
\end{equation}
for all $\alpha\in\Delta$.
We note that \eqref{weakrefs} is quite different from
\begin{equation}\label{strongrefs}
S_\alpha-2S_\alpha\subset S_{\alpha}.
\end{equation}
We call a subset satisfying \eqref{weakrefs}
a {\bf reflection space} of $G$,
and satisfying
\eqref{strongrefs}
a {\bf symmetric reflection space} of $G$
(though we used \eqref{strongrefs} to be the definition
of a {\bf refection space} of $G$ in \cite{Y2} and \cite{NY}).
Also, we call {\bf full} if such a subset generates $G$.
A reflection space $E$ is called {\bf pointed} if $0\in E$.

We will not classify the family $\{S_\alpha\mid\alpha\in\Delta\}$
of supports,
which is maybe difficult to classify.
We only concentrate on each reflection space $S_\alpha$.

\example
If $G=\Bbb Z$,
then $p\Bbb Z+e$ 
for any $p, e\in\Bbb Z$ 
is a reflection space,
and it is full if and only if $(p,e)=1$.
Note that any singleton $\{e\}$ is a reflection space.
On the other hand,
a symmetric reflection space of $\Bbb Z$ is just $p\Bbb Z$ or $p(2\Bbb Z+1)$,
and $2\Bbb Z+1$ and $\Bbb Z$ are the only full symmetric reflection spaces of $\Bbb Z$.
(We will summarize these in Proposition \ref{1ref}.)
\endexample

\example\label{linearsol}
Let $A$ be a matrix.
Then the solution space to the system $A{\bf x}={\bf b}$ of equations 
is a reflection space.
In fact, if $A{\bf x}={\bf b}$ and $A{\bf y}={\bf b}$,
then $A(2{\bf x}-{\bf y})=2A{\bf x}-A{\bf y}=2{\bf b}-{\bf b}={\bf b}$.
Hence $2{\bf x}-{\bf y}$ is also a solution.
\endexample

\medskip

The following lemmas are basic.

\lemma\label{g+Eref}
If $E$ is a reflection space of $G$, then 
$-E$ and $E+g$ are also reflection spaces
for any $g\in G$.
In particular, if $e\in E$,
then $E-e$ for any $e\in E$ is a pointed reflection space.

\endlemma
\proof
For $x,\  y\in E$,
we have $2(-x)-(-y)=-(2x-y)\in -E$.
We also have $2(x+g)-(y+g)=2x-y+g\in E+g$.
\qed

\lemma\label{basic0}
Let $E$ be a symmetric reflection space of $G$.
Then 
$-E=E$.
Hence $E+2E\subset E$ and $2E-E\subset E$. 
Thus a symmetric reflection space is a reflection space.
\endlemma
\proof
For $x\in E$,
we have $x-2x=-x\in E$.
Hence $-E\subset E$.
Thus $E\subset -E$.
\qed

\lemma\label{basic}
Let $E$ be a reflection space of $G$.
Then 
$$
\text{$E$ is pointed $\Longrightarrow$ $E$ is a symmetric reflection space.}
$$
Hence, a pointed reflection space is a pointed symmetric reflection space.
\endlemma
\proof
Since $0\in E$, we get $-E\subset E$.
Hence $E-2E=-(2E-E)\subset E$.
\qed

\lemma\label{uniqueness}
Let $E$ be a subset of $G$.
For any $e,e'\in E$,
we have
$$\langle E-e\rangle =\langle E-e'\rangle,$$
where the bracket
$\langle A\rangle$
denotes the subgroup generated by a subset $A$ of $G$.
\endlemma
\proof
For $x\in E$,
we have $x-e, e'-e\in E-e$.
Hence $x-e'=x-e-(e'-e)\in \langle E-e\rangle$,
and so
$\langle E-e'\rangle\subset \langle E-e\rangle$.
Similarly, we have
$\langle E-e\rangle\subset \langle E-e'\rangle$.
\qed

\lemma\label{containcyc}
Let $E$ be a pointed reflection space of $G$,
and let $e\in E$.
Then $\langle e\rangle\subset E$.
\endlemma
\proof
Since $0\in E$
(so $E$ is symmetric),
we have $\pm 2e=0\pm 2e\in E$ 
and $\pm 3e=\pm (e+2e)\in E$.
Similarly,
we have $2me=0\pm (2e+\cdots +2e)\in E$ and 
$(2m+1)e=e\pm (2e+\cdots +2e)\in E$
for all $m\in\Bbb Z$.
\qed

More generally, we have:
\lemma\label{containdouble}
Let $E$ be a symmetric reflection space of $G$.
Suppose that $\{e_i\}_{i\in\frak I}\subset E$,
where $\frak I$ is any index set.
Then $E+2\langle e_i\rangle_{i\in\frak I}\subset E$.
Hence,
$E+2\langle E\rangle= E$.

\endlemma
\proof
Let $x\in E+2\langle e_i\rangle_{i\in\frak I}$.
Then $x=e+2\sum_{j=1}^n \epsilon_je_{i_j}$,
where
$\epsilon_j=1$ or $-1$,
and $e_{i_j}\in\{e_i\}_{i\in\frak I}$.
Thus
$x=e+2\epsilon_1e_{i_1}+\cdots+2\epsilon_ne_{i_n}\in E$, inductively.
(Note that $-e_{i_j}\in E$ by
Lemma \ref{basic0}).
\qed

\enskip

Now we classify reflection spaces of $\Bbb Z$.

\proposition\label{1ref}
Let $E$ be a subset of $\Bbb Z$.
Then
\begin{equation}\label{prefsc}
\text{$E$ is a pointed reflection space
$\Longrightarrow$
$E=p\Bbb Z$}
\end{equation}
for some $p\in\Bbb Z_{\geq 0}$.
So a pointed reflection space of $\Bbb Z$ is just a subgroup of $\Bbb Z$.

Moreover,
\begin{equation}\label{weakrefsc}
\text{$E$ is a reflection space
$\Longrightarrow$
$E=p\Bbb Z+e$}
\end{equation}
for any $e\in E$ and some $p\in\Bbb Z_{\geq 0}$,
and
\begin{equation}\label{refsc}
\text{$E$ is a symmetric reflection space
$\Longrightarrow$
$E=p\Bbb Z$ or $p(2\Bbb Z+1)$}.
\end{equation}
\endproposition

\proof
Since $\langle E\rangle$ is a subgroup of $\Bbb Z$,
we have $\langle E\rangle=p\Bbb Z$ for some $p\in\Bbb Z_{\geq 0}$.
Thus it is enough to show that $p\in E$,
by Lemma \ref{containcyc}.
Let $p=\sum_i e_i$ for $e_i\in E\subset p\Bbb Z$.
If all $e_i= 2pk_i$ for $k_i\in\Bbb Z$,
then
$p=2p\sum_i k_i$,
and so $1=2\sum_i k_i$,
which is absurd.
Hence for some $j$,
we have $e_j=p(2k_j+1)\in E$.
Since $pk_j\in\langle E\rangle$,
we have,
by Lemma \ref{containdouble},
$p=p(2k_j+1)-2pk_j\in E$.

For \eqref{weakrefsc},
note that
$E-e$ for $e\in E$ is a pointed reflection space $\Bbb Z$,
by Lemma \ref{g+Eref}.
Hence by \eqref{prefsc}, we have
$E-e=p\Bbb Z$ for some $p\in\Bbb Z_{\geq 0}$.
Thus $E=p\Bbb Z+e$.

For \eqref{refsc},
we have $E=p\Bbb Z+e$, by the same reason above.
But if $p=0$, then $e$ has to be $0$,
and if $p=1$, then $E=\Bbb Z$.
Thus, we may assume $p>1$
and also $e>0$.
Let $d=(p,e)$ (the greatest common divisor).
Then
one can write $E=d(p'\Bbb Z+e')$ so that $(p',e')=1$
for some $p'>1$.
Since $-E\subset E$,
we have $de'\equiv -de'$ (mod $dp'=p$).
Hence $p'\mid 2e'$,
and
we obtain $p'=2$.
Then one can take $e'$ to be $1$.
Therefore, $E=d(2\Bbb Z+1)$.
\qed

\medskip

In general, we have the following.

\proposition\label{refingeneral}
Let $E$ be a subset of $G$.
Then
\begin{equation}\label{prefsg}
\text{$E$ is a symmetric reflection space
$\Longrightarrow$
$E=\bigcup_{i=1}^m\  (2\langle E\rangle+e_i)$}
\end{equation}
for some $1\leq m\leq |\langle E\rangle/2\langle E\rangle|$
(possibly infinite)
and some
$e_i\in E$,
and if $E$ is pointed,
then some $e_i=0$.

Moreover,
\begin{equation}\label{weakrefsg}
\text{$E$ is a reflection space
$\Longrightarrow$
$E=\bigcup_{i=1}^m\  (2\langle E-e\rangle+e_i)$}
\end{equation}
for some $1\leq m\leq \infty$ and any $e\in E$ (see Lemma \ref{uniqueness}),
and some $e_i\in E$
(possibly $e_i\notin \langle E-e\rangle$).

Conversely,

{\rm (i)}\ \ 
$E=\bigcup_{i=1}^m\  (2S+s_i)$
for 
any subgroup $S$ of $G$,
 $s_i\in S$, 
 and any $1\leq m\leq |S/2S|$
 (possibly infinite)
is a symmetric reflection space.

{\rm (ii)}\ \
Let $E=\bigcup_{i=1}^m\  (S+x_i)$
for any subgroup $S$ of $G$
and some $x_i\in G$.
Suppose also that
for all $1\leq i, j\leq m$  $(1\leq i, j<\infty$ if $m=\infty)$,
there exists some $1\leq k\leq m$
 $(1\leq k<\infty$ if $m=\infty)$
such that
$S+2x_i-x_j=S+x_k$.
Then $E$ is a reflection space.
(If $m=1$, then $E$ is always a reflection space by Lemma \ref{g+Eref}.)

\endproposition
\proof
For \eqref{prefsg},
it follows from
Lemma \ref{containdouble}. 

For \eqref{weakrefsg}, since
$E-e$ is a (pointed) symmetric reflection space,
we have, by \eqref{prefsg},
$$
E-e=\bigcup_{i=1}^m\  (2\langle E-e\rangle+g_i)
$$
for some $g_i\in E-e$.
So letting $e_i:=g_i+e\in E$,
we obtain \eqref{weakrefsg}.

Conversely,
for (i),
let
$2s+s_i,\  2s'+s_j\in E$.
Then
$2s+s_i-2(2s'+s_j)\in 2S+s_i$.
Hence $E$ is a symmetric reflection space.

For (ii),
let
$s+x_i,\  s'+x_j\in E$.
Then we have
$$2(s+x_i)-(s'+x_j)
=2s-s'+2x_i-x_j\in S+x_k,
$$
by our assumption.
Hence $E$ is a reflection space.
\qed

\medskip

\definition
Let $E$ and $E'$ be reflection spaces of abelian groups $G$ and $G'$, respectively.
We say that $E$ is {\bf isomorphic} to $E'$, denoted $E\cong E'$,
if there exists a bijection $f:E\longrightarrow E'$
such that $f(2x-y)=2f(x)-f(y)$ for all $x,y\in E$.
In particular,
if there exists a group isomorphism $\varphi:G\longrightarrow G'$
such that $\varphi(E)=E'$,
then $E\cong E'$.

\enddefinition

\lemma
Let $E$ be a reflection space of an abelian group $G$.
Then $E+g\cong E$ for any $g\in G$.
\endlemma
\proof
Let $f:E\longrightarrow E+g$ be the bijection
defined by $f(x)=x+g$ for $x\in E$.
Then we have $2f(x)-f(y)=2(x+g)-(y+g)=2x-y+g=f(2x-y)$,
and hence $f$ is an isomorphism.
\qed

(Note that if $g\neq 0$,
the translation $f$ cannot be the restriction of an automorphism of $G$
since $f(0)\neq 0$.)

\example\label{exref}
(1)
For any $(p_1, p_2), (e_1,e_2)\in\Bbb Z^2$,
$$E= (p_1\Bbb Z+e_1)\times(p_2\Bbb Z+e_2)$$
is a reflection space of $\Bbb Z^2$,
which is isomorphic to $p_1\Bbb Z\times p_2\Bbb Z$.

(2)
Taking $S=3\Bbb Z^2$ in Proposition \ref{refingeneral}(ii),
\begin{align*}
E_1
&=
\big(3\Bbb Z^2+(1,0)\big)
\cup \big(3\Bbb Z^2+(0,1)\big)
\cup \big(3\Bbb Z^2+(2,2)\big)\quad\text{or}\\
E_2&=
\big(3\Bbb Z^2+(2,0)\big)
\cup \big(3\Bbb Z^2+(0,1)\big)
\cup \big(3\Bbb Z^2+(1,2)\big)
\end{align*}
is a reflection space of $\Bbb Z^2$.

(3)
Taking $S=6\Bbb Z^2$ in Proposition \ref{refingeneral}(ii),
\begin{align*}
E
&=
\big(6\Bbb Z^2+(1,0)\big)
\cup \big(6\Bbb Z^2+(0,1)\big)
\cup \big(6\Bbb Z^2+(2,5)\\
&\cup
\big(6\Bbb Z^2+(5,2)\big)
\cup \big(6\Bbb Z^2+(3,4)\big)
\cup \big(6\Bbb Z^2+(4,3)\big)
\end{align*}
is a reflection space of $\Bbb Z^2$.

\endexample

\definition\label{refgensub}
For a subset $A\subset G$,
the reflection space generated by $A$
is denoted by
$ [A]$,
and 
the symmetric reflection space generated by $A$
is denoted by
$[A\rangle$.
Also,
the pointed reflection space generated by $A$
is denoted by
$[A\rangle_0$.

More precisely, we have
$$
[A]=\{g_1\cdot(g_2\cdots\cdot (g_r\cdot g_{r+1}))\cdot\cdot)\mid g_i\in A\},
$$
$$
[A\rangle=\{g_1\circ(g_2\cdots\circ (g_r\circ g_{r+1}))\cdot\cdot)\mid g_i\in A\},
$$
$$
\text{and}\quad [A\rangle_0=[A\cup\{0\}]\ (=[A\cup\{0\}\rangle\  \text{by Lemma \ref{basic}}),
$$
where
$g_i\cdot g_j=2g_i-g_j$ and $g_i\circ g_j=g_i-2g_j$.

\enddefinition

\example\label{}
(1)
Let $E=[6,15]$ be the reflection space generated by $6$ and $15$ in $\Bbb Z$.
Since $E=p\Bbb Z+e$ for some $p,e\in\Bbb Z$,
we have
$15-6=9=pm$
for some $m\in\Bbb Z$.
So $p$ can be $1$, $3$, or $9$.
But $E$ should be the smallest one containing $6$ and $15$,
and hence $E=9\Bbb Z+6$.
Similarly, for $a,b\in\Bbb Z$,
one can show that
\begin{equation}\label{formula1}
[a,b]=(b-a)\Bbb Z+a.
\end{equation}

(2)
Let $S=[6,15\rangle$ be the symmetric reflection space generated by $6$ and $15$ in $\Bbb Z$.
Since $\langle S\rangle=3\Bbb Z$,
we have
$S=3\Bbb Z$ or $6\Bbb Z+3$.
But $6\notin 6\Bbb Z+3$, we get $S=3\Bbb Z$.

(3)
Let $S=[15,27\rangle$ be the symmetric reflection space generated by $15$ and $27$ in $\Bbb Z$.
Since $\langle S\rangle=3\Bbb Z$,
we have
$S=3\Bbb Z$ or $6\Bbb Z+3$.
Since $S$ should be the smallest one,
we get $S=6\Bbb Z+3$.
Note that $[6,15\rangle_0=3\Bbb Z$.

\endexample

We generalize the formula \eqref{formula1}.
\begin{proposition}
Let $x,y\in G$ for an abelian group $G$.
Then we have 
\begin{equation}\label{genformula0}
[x,y]=\langle y-x\rangle+x
\end{equation}
and so $[x,y]-x=\langle x-y\rangle$.
Also, we have
\begin{equation}\label{genformula00}
[x,y]=\langle x-y\rangle+x=\langle x-y\rangle+y.
\end{equation}

\end{proposition}
\proof
The right-hand side is a reflection space containing $x$ and $y$
(see Lemma \ref{g+Eref}),
and so it is enough to show that 
$[x,y]\supset\langle y-x\rangle+x$.
We show that
$m(y-x)+x\in [x,y]$ for all $m\in\Bbb Z$.
It is clear that $m=0, \pm 1$.
Assume that $|m|>1$ and we use the induction on $m$.
If $m$ is even,
then $|m/2|<|m|$,
and so we have
$\frac{m}{2}(y-x)+x\in [x,y]$.
Hence $m(y-x)+2x-x=m(y-x)+x\in [x,y]$.
If $m$ is odd,
then $|\frac{m+1}{2}|<|m|$,
and so we have
$\frac{m+1}{2}(y-x)+x\in [x,y]$.
Hence $(m+1)(y-x)+2x-y=m(y-x)+y-x+2x-y
=m(y-x)+x\in [x,y]$.
Therefore, \eqref{genformula0} holds,
and
\eqref{genformula00} is clear since $[x,y]=[y,x]$.
\qed

\example\label{twogenerators}
Let $L=\frak g\otimes F[t_1^{\pm 1},\ldots, t_n^{\pm 1}]$ be the multi-loop algebra,
where $\frak g$ is a finite-dimensional split simple Lie algebra.
Let $e$ and $f$ be nonzero root vectors of $\frak g$
for roots $\alpha$ and $-\alpha$, respectively.
Let
\begin{equation*}\label{generators}
A:=\{x,y\}\subset \Bbb Z^n\quad\text{and}\quad
U:=
\{e\otimes t^{x},\quad
f\otimes t^{-x},\quad
e\otimes t^{y},\quad
f\otimes t^{-y}
\},
\end{equation*}
where
$t^{v}$ means $t_1^{v_1}\cdots t_n^{v_n}$
for $v=(v_1,\ldots, v_n)\in\Bbb Z^n$.

Let $M$ be the subalgebra of $L$
generated by $U$.
Then
$M$ is a general Lie torus of type $\text A_1=\{\pm\alpha\}$.
Since $\supp_{\Bbb Z^n} M_\alpha$
is a reflection space containing $A$,
we have $[A]\subset \supp_{\Bbb Z^n} M_\alpha$.
For the other inclusion,
let $z\in\supp_{\Bbb Z^n} M_\alpha$.
Since $M$ is generated by $U$,
one can see that
$z=(m+1)x-mx=x$, $(m+1)x-my$, $(m+1)y-mx$ or $(m+1)y-my=y$
for some $m\in\Bbb Z_{\geq 0}$.
But each of them is in $\Bbb Z(x-y)+x=[A]$.
Hence we have
$\supp_{\Bbb Z^n} M_\alpha=[A]$.
Similarly,
we have $\supp_{\Bbb Z^n} M_{-\alpha}=-[A]$.
Note that
$\supp_{\Bbb Z^n} M_0=\supp_{\Bbb Z^n} M_\alpha+\supp_{\Bbb Z^n} M_{-\alpha}=[A]-[A]=\Bbb Z(x-y)$.

\medskip

Let $s\in\hom(\langle\alpha\rangle,\Bbb Z^n)$,
defined by
$s(\alpha)=x$.
Then the $s$-isotope $M^{(s)}$ 
is a normal Lie $1$-torus.
In fact, 
we have
$$
Fe\otimes t^{x}=M_{\alpha}^{x}
=(M^{(s)})_{\alpha}^{0}\quad\text{and}
\quad
Ff\otimes t^{-x}
=M_{-\alpha}^{-x}
=(M^{(s)})_{-\alpha}^{0},
$$
and letting $p:=y-x$,
$$
Fe\otimes t^{y}=M_{\alpha}^{y}
=M_{\alpha}^{p+x}
=(M^{(s)})_{\alpha}^{p}\quad\text{and}
\quad
Ff\otimes t^{-y}=M_{-\alpha}^{-y}
=M_{-\alpha}^{-p-x}
=(M^{(s)})_{-\alpha}^{-p},
$$
and so $M^{(s)}$ is isograded isomorphic
to the loop algebra $\sll_2(F[X^{\pm 1}])$,
where $X=t^p$.

\endexample

\begin{lemma}
Let $x,y\in G$ for an abelian group $G$.
Then we have 
\begin{equation}\label{genformula}
[x,y\rangle=(2\langle x,y\rangle+x)\cup(2\langle x,y\rangle+y)
\end{equation}
\begin{equation}\label{genformula2}
\text{and}\quad
[x,y\rangle_0=2\langle x,y\rangle\cup [x,y\rangle.
\end{equation}

\end{lemma}
\proof
It follows from Proposition \ref{refingeneral}.
\qed

\lemma
We have the formula
\begin{equation}\label{niceformula}
[x,y\rangle
=\langle x-y\rangle+\langle x+y\rangle+x
=\langle x-y\rangle+\langle x+y\rangle+y
=[x,y]+\langle x+y\rangle.
\end{equation}
\endlemma
\proof
By \eqref{genformula00},
it is enough to show the first identity.
From \eqref{genformula}
and the inclusion
$$
2\langle x,y\rangle\subset\langle x-y\rangle+\langle x+y\rangle
$$
(since $2mx+2ny=(m-n)(x-y)+(m+n)(x+y)$ for $m,n\in\Bbb Z$),
we have $[x,y\rangle\subset\langle x-y\rangle+\langle x+y\rangle+x$.
We show the other inclusion.
For $X:=m(x-y)+n(x+y)+x$,
if $m+n$ is even,
then $-m+n$ is also even,
and so $X=(m+n)x+(-m+n)y+x\in[x,y\rangle$.
If $m+n$ is odd,
then $m+n+1$ and $-m+n-1$ are even,
and so
$$X=(m+n+1)x+(-m+n-1)y+y\in[x,y\rangle.$$
\qed

\example\label{pmtwogenerators}
In the notations in Example \ref{twogenerators},
Let 
\begin{equation*}\label{pmgenerators}
T:=U\cup
\{e\otimes t^{-x},\quad
f\otimes t^{x},\quad
e\otimes t^{-y},\quad
f\otimes t^{y}
\}.
\end{equation*}

Let $N$ be the subalgebra of $L$
generated by $T$.
Then
$N$ is a general Lie torus of type $\text A_1=\{\pm\alpha\}$.
Since $\supp_{\Bbb Z^n} N_\alpha$
is a reflection space containing $\pm \ A$,
we have $[A\rangle\subset \supp_{\Bbb Z^n} N_\alpha$.
For the other inclusion,
let $z\in\supp_{\Bbb Z^n} N_\alpha$.
Since $N$ is generated by $T$,
we have
$z\in \epsilon_1+\cdots +\epsilon_{2m+1}$
for some $m\in\Bbb Z_{\geq 0}$,
where all
$\epsilon_i=\epsilon=\{\pm x, \pm y\}$.
We show that 
$\epsilon_1+\cdots +\epsilon_{2m+1}\subset [A\rangle$,
by induction on $m$.
It is clear for $m=0$.
Suppose that the statement is true for $m$.
Let $p\in\epsilon_1+\cdots +\epsilon_{2m+1}$.
Then we have
$$p+\epsilon+\epsilon=\{p, p+2x,p-2x, p+2y, p-2y, p+x+y, p+x-y, p-x+y, p-x-y\}.$$
The first 5 elements in the right hand side are clearly in $[A\rangle$.
The last 4 elements are also in $[A\rangle$, by \eqref{niceformula}.
Hence we have shown
$\supp_{\Bbb Z^n} N_\alpha=[A\rangle$.
Similarly,
we have $\supp_{\Bbb Z^n} N_{-\alpha}=[A\rangle$.
Note that
$\supp_{\Bbb Z^n} N_{0}=\supp_{\Bbb Z^n} N_\alpha+\supp_{\Bbb Z^n} N_{-\alpha}
=[A\rangle+[A\rangle=\langle x+y,\ x-y\rangle\neq \langle A\rangle$.

\medskip

Let $s\in\hom(\langle\alpha\rangle,\Bbb Z^n)$,
defined by
$s(\alpha)=x$.
Then the $s$-isotope $N^{(s)}$ 
is a normal Lie $2$-torus.
In fact, 
we have
$$
Fe\otimes t^{x}=N_{\alpha}^{x}
=(N^{(s)})_{\alpha}^{0}\quad\text{and}
\quad
Ff\otimes t^{-x}
=N_{-\alpha}^{-x}
=(N^{(s)})_{-\alpha}^{0},
$$
and letting $p:=y-x$
and $q:=x+y$,
$$
Fe\otimes t^{y}=N_{\alpha}^{y}
=N_{\alpha}^{p+x}
=(N^{(s)})_{\alpha}^{p}\quad\text{and}
\quad
Ff\otimes t^{-y}=N_{-\alpha}^{-y}
=N_{-\alpha}^{-p-x}
=(N^{(s)})_{-\alpha}^{-p},
$$
$$
Fe\otimes t^{-y}=N_{\alpha}^{-y}
=N_{\alpha}^{-q+x}
=(N^{(s)})_{\alpha}^{-q}\quad\text{and}
\quad
Ff\otimes t^{y}=N_{-\alpha}^{y}
=N_{-\alpha}^{q-x}
=(N^{(s)})_{-\alpha}^{q},
$$
$$
Fe\otimes t^{-x}=N_{\alpha}^{-x}
=N_{\alpha}^{p-q+x}
=(N^{(s)})_{\alpha}^{p-q}\quad\text{and}
\quad
Ff\otimes t^{x}=N_{-\alpha}^{x}
=N_{-\alpha}^{-p+q-x}
=(N^{(s)})_{-\alpha}^{q-p}.
$$
So
we have $\supp_{\Bbb Z^n}N^{(s)}_\alpha=[0, p, - q, p-q]=\langle p, q\rangle=[0,- p, q, q-p]
=\supp_{\Bbb Z^n}N^{(s)}_{-\alpha}$.
Thus $N^{(s)}$ is isograded isomorphic
to the double loop algebra $\sll_2(F[X^{\pm 1},Y^{\pm 1}])$,
where $X=t^p$ and $Y=t^q$.
\endexample

\remark
We note that
$[A,0]=[ A\rangle_0$
(see Lemma \ref{basic}).
However, for
the subalgebra $N'$ of $L$ generated by
$U\cup\{e, f\}$,
which is a normal Lie $2$-torus,
we have $\supp_ {\Bbb Z^n}N'_\alpha\neq [ A\rangle_0$.
In fact, $\supp_ {\Bbb Z^n}N'_\alpha$ contains
$\epsilon_1+\cdots +\epsilon_{m}$
for any $m\in\Bbb Z_{\geq 0}$,
where $\epsilon_i$ is defined above.
Hence
$\supp_ {\Bbb Z^n}N'_\alpha$ contains $ \langle A\rangle$,
and
therefore,
$\supp_ {\Bbb Z^n}N'_\alpha= \langle A\rangle$.

There are examples that a normal Lie torus $W$ of type $\text A_1$
satisfies $\supp_ {\Bbb Z^n}W_\alpha=[ A\rangle_0$.
For example, let
$W=\text{TKK}\big(H(F_h[t_1^{\pm 1},t_2^{\pm 1}],*)\big)$
be the Tit-Koecher-Kantor Lie algebra
constructed from 
$H(F_h[t_1^{\pm 1},t_2^{\pm 1}],*)$,
which is the Jordan algebra of the fixed points
in the quantum torus
$F_h[t_1^{\pm 1},t_2^{\pm 1}]$ defined by $t_2t_1=-t_1t_2$,
by the involution $*$ determined by $t_1^*=t_1$ and $t_2^*=t_2$.
Let $x=(1,0)$, $y=(0,1)$ and $A=\{x, y\}$.
Then
$\supp_ {\Bbb Z^2}W_\alpha=[ A\rangle_0$.
\endremark

\enskip

\enskip

\section{A new definition of a Lie torus}
\label{}

Let us recall the definition of a locally extended affine root system.
\definition
Let $V$ be a vector space over $\Bbb Q$
with a positive semidefinite symmetric bilinear form $(\cdot,\cdot)$. 
A subset $\mathfrak R$ of $V$ 
is called a {\it locally extended affine root system}
or a {\it LEARS} for short if
$\frak R$ satisfies the following:

\item[(A1)]  $(\alpha,\alpha)\neq 0$  for all $\alpha\in\mathfrak R$, and
$\mathfrak R$ spans $V$;

\item[(A2)]
$\langle\alpha,\beta\rangle\in\Bbb Z$ for all 
$\alpha, \beta\in\mathfrak R$,
where $\langle\alpha,\beta\rangle=\frac{2(\alpha,\beta)}{(\beta,\beta)}$;

\item[(A3)]
$\sigma_\alpha(\beta)\in \mathfrak R$ for all $\alpha,\beta\in \mathfrak R$, where
$\sigma_\alpha(\beta)=\beta-\langle\beta,\alpha\rangle\alpha$;

\item[(A4)]
$\mathfrak R=\mathfrak R_1\cup \mathfrak R_2$ and $(\mathfrak R_1,\mathfrak R_2)=0$ imply $\mathfrak R_1=\emptyset$ or
$\mathfrak R_2=\emptyset$\ \  \
(irreducibility).

\medskip

A LEARS $\mathfrak R$ is called {\it reduced} if
$2\alpha\notin \mathfrak R$ for all $\alpha\in\mathfrak R$.

\enddefinition

Note that if $V$ is finite-dimensional and $(\cdot,\cdot)$ is positive definite,
then $\mathfrak R$ is exactly a finite irreducible root system
(see [MY1, Prop. 4.2]).

\medskip

Now we define a new Lie torus determined by
a LEARS.

\definition\label{defnewL}
Let $(\frak R, V)$ be a LEARS,
and let 
$Q:=\langle\frak R\rangle$ be the root lattice,
i.e.,
the subgroup of $V$
generated by $\frak R$.
Let
$$L=\bigoplus_{\xi\in Q}\ L_\xi$$
be
a $Q$-graded Lie algebra over $F$ generated by
$$\bigcup_{\alpha\in\frak R} \ L_\alpha
\quad\text{so that}\quad
(\supp L)^\times=\frak R,$$
where
$(\supp L)^\times=\{\xi\in\supp L\mid (\xi,\xi)\neq 0\}$.

(1)\ \
$L$
is called {\bf division graded}
if for each $\alpha\in\frak R$ and $0\neq x\in L_\alpha$,
there exists $y\in L_{-\alpha}$
such that
$$\alpha^\vee : = [x, y] \in L_0
\quad\text{satisfies}\quad
[\alpha^\vee , z] = \langle \xi,\alpha\rangle z$$ for all  $z \in L_\xi$
for $\xi \in Q$.

(2)\
A division graded Lie algebra $L=\oplus_{\xi\in Q}\ L_\xi$
is called a {\bf Lie $\frak R$-torus} if
$$\dim_F L_\alpha=1$$ for all $\alpha\in\frak R$.

\enddefinition

\lemma
A Lie $\frak R$-torus $L$ only exists for a reduced LEARS $\frak R$.

\endlemma

\proof
Suppose $\alpha,2\alpha\in\frak R$.
Then we have $\dim L_\alpha=\dim L_{2\alpha}=1$.
Let
$A:=L_{-\alpha}\oplus L_{0}\oplus L_{\alpha}\cong\sll_2(F)$,
and
$M:=L_{-2\alpha}\oplus L_{-\alpha}\oplus L_{0}\oplus L_{\alpha}
\oplus L_{2\alpha}$
the $A$-submodule of $L$.
Then by the complete reducibility of $M$,
we get $\dim L_\alpha=2$
since $M$ decomposes the direct sum of two irreducible submodules.
This
 is a contradiction.
\qed

\enskip

Let $$V^0:=\{x\in V\mid (x,y)=0\ \text{for all $y\in V$}\}$$ be the radical of the form.
Note that
$$V^0=\{x\in V\mid (x,x)=0\}.$$
We call $\dim_\Bbb Q V^0$ the {\it null dimension} of $\mathfrak R$,
which can be any cardinality.
We call a LEARS
$(\mathfrak R,V)$ an {\it extended affine root system}
or an {\it EARS} for short
if 
$$\text{$\dim_\Bbb QV/V^0<\infty$ and
$\langle \mathfrak R\rangle$ is free.}
$$
This coincides with the concept, which was first introduced by Saito in 1985 \cite{S}.
The notion of an EARS was also used in a different sense in \cite{AABGP},
but Azam showed that there is a natural correspondence between the two notions in \cite{A}.
We use here the Saito's one since he is the first person who defined it and his root systems naturally 
generalize Macdonald's {\it affine root systems} in \cite{M}.

\medskip

Let $(\mathfrak R,V)$ be a LEARS, and $(\bar{ \mathfrak R},\bar V)$ the canonical image onto $V/V^0$. 
Then $\bar V$ admits the induced positive definite form, and  thus
$$
\text{$(\bar {\mathfrak R},\bar V)$ is a locally finite irreducible root system.}
$$

\medskip

Now we show the new Lie $\mathfrak R$-torus 
$L=\bigoplus_{\xi\in Q}\ L_\xi$
is a general Lie torus.
Let $V'$ be a section of $\bar V$, i.e.,
$V=V'\oplus V^0$,
and let
$$
\Delta:=\{\dot\alpha\in V'\mid \bar {\dot\alpha}\in\bar{\frak R}\}.
$$
Then $(\Delta, V')$ is a locally finite irreducible root system
isomorphic to $\bar{\frak R}$.
For $\dot\alpha\in\Delta$, let
$$
S_{\dot\alpha}:=\{v\in V^0\mid \dot\alpha+v\in\frak R\},
$$
and
let
$$G:=\bigg\langle\bigcup_{\dot\alpha\in\Delta}\ S_{\dot\alpha}\bigg\rangle.$$
So for $\xi\in Q$, we have
$$\xi=\sum_{\alpha\in\frak R} a_\alpha\alpha
=\sum_{\alpha\in\frak R} a_\alpha(\dot\alpha+g_{\alpha})
=\sum_{\dot\alpha\in\Delta} b_{\dot\alpha}\dot\alpha
+\sum_{\dot\alpha\in\Delta} k_{\dot\alpha}
$$
for $a_\alpha\in\Bbb Z$, \
 $g_{\alpha}\in S_{\dot\alpha}$,\ \
$b_{\dot\alpha}=\sum_{\alpha\in\bar{\dot\alpha}}a_\alpha$
and $k_{\dot\alpha}=\sum_{\alpha\in\bar{\dot\alpha}}a_\alpha g_\alpha$.
Note that
$$\text{
$(\xi,\xi)\neq 0$ $\Leftrightarrow$ the first term
$\sum_{\dot\alpha\in\Delta} b_{\dot\alpha}\dot\alpha\neq 0$.}
$$
Note also that
$$\sum_{\dot\alpha\in\Delta} b_{\dot\alpha}\dot\alpha\in\Delta\
\Leftrightarrow\
\overline{\sum_{\dot\alpha\in\Delta} b_{\dot\alpha}\dot\alpha}\in\bar{\frak R}$$
since
$\sum_{\dot\alpha\in\Delta} b_{\dot\alpha}\dot\alpha\in V'$.
Therefore,
since $(\supp L)^\times=\frak R$,
we have
$$\xi\in\supp L\
\Leftrightarrow\
\sum_{\dot\alpha\in\Delta} b_{\dot\alpha}\dot\alpha\in\Delta\cup\{0\}.$$
Let 
$\mu:=\sum_{\dot\alpha\in\Delta} b_{\dot\alpha}\dot\alpha$
and $g:=\sum_{\dot\alpha\in\Delta} k_{\dot\alpha}$,
which are unique for $\xi$.
Then through
$$\xi\in\supp L\ \longrightarrow\  (\mu, g)\in (\Delta\cup\{0\})\times G,$$
we have
$$L=\bigoplus_{\xi\in Q}\ L_\xi
=\bigoplus_{(\mu, g)\in (\Delta\cup\{0\})\times G}\ L_\mu^g,$$
where
$L_\mu^g=L_{\mu+g}$ if $\mu+g\in Q$
and $L_\mu^g=0$ otherwise.
Then it is clear that $L$ satisfies the axioms of a general Lie $G$-torus,
except (LT2).

We show (LT2).
For every $g \in G$, we need to show that 
$L_0^g \subset  \sum_{\mu \in \Delta,\ h \in G} \  [L_\mu^h, L_{-\mu}^{g-h}]$.
(The other inclusion is clear.)
Since $L_0^g=L_{0+g}$ is contained in the subalgebra generated by
$\bigcup_{\alpha\in\frak R} \ L_\alpha$
and in the degree $(0+g)$-space,
we have
$$
L_0^g 
\subset\bigg(\sum_{\alpha,\beta\in\frak R}\ [L_\alpha, L_\beta]\bigg)_{0+g}
=  \sum_{\dot\alpha \in \frak R, \ g_\alpha\in G} \  [L_{\dot\alpha+g_{\alpha}}, L_{-\dot\alpha+g-g_{\alpha}}]
=  \sum_{\mu \in \Delta,\ h \in G} \  [L_\mu^h, L_{-\mu}^{g-h}].
$$

Therefore:
\begin{theorem}
A Lie $\frak R$-torus can be identified with a general Lie $G$-torus,
where $G$ is determined by a section $V'$ above.
\end{theorem}

Let $\Pi$ be a reflectable base of $\bar{\frak R}$.
For each $\alpha\in\Pi$,
let $\dot\alpha$ be an element of $\frak R$ so that $\bar{\dot\alpha}=\alpha$.
Let
$$
V':=\spa\{\dot\alpha\in\frak R\mid \alpha\in\Pi\}.
$$
Since $\Pi$ is a basis of $\bar V$,
we have $V=V'\oplus V^0$.
We call this $V'$ a {\bf reflectable section} relative to a reflectable base $\Pi$
(and a choice of $\{\dot\alpha\in\frak R\mid \alpha\in\Pi\}$).
We set $\dot\alpha\in V'$ for other $\alpha\in\bar{\frak R}$
so that $\bar{\dot\alpha}=\alpha$.

\claim\label{refsection}
If $\bar{\dot\alpha}\in\bar{\frak R}^\red$,
then $\dot\alpha\in\frak R$.
Hence $0\in S_{\dot\alpha}$.
\endclaim
\proof
We have
$\bar{\dot\alpha}=\sigma_{\alpha_1}\cdots\sigma_{\alpha_k}(\alpha_{k+1})$
for some $\alpha_1, \ldots, \alpha_{k+1}\in\Pi$.
Then we have
$$
\sigma_{\dot\alpha_1}\cdots\sigma_{\dot\alpha_k}(\dot\alpha_{k+1})\in V'
\quad\text{and}\quad
\overline{\sigma_{\dot\alpha_1}\cdots\sigma_{\dot\alpha_k}(\dot\alpha_{k+1})}
=\bar{\dot\alpha}.
$$
Hence
$\dot\alpha=\sigma_{\dot\alpha_1}\cdots\sigma_{\dot\alpha_k}(\dot\alpha_{k+1})
\in\frak R$
since each $\dot\alpha_i\in\frak R$.
\qed

\medskip

Thus
$L$ satisfies (LT6), and hence:
\begin{theorem}
A Lie $\frak R$-torus can be identified with a normal Lie $G$-torus,
where $G$ is determined by a reflectable section.
\end{theorem}

\medskip

Suppose that $G$ is a torsion-free abelian group.
Then $\langle\Delta\rangle\times G$ is also a torsion-free abelian group
for a locally Lie $G$-torus
$L=\bigoplus_{\mu\in\Delta\cup\{0\}}\ \bigoplus_{g\in G}\ L_{\mu}^g$.
Hence $\langle\Delta\rangle\times G$ can embed into the vector space 
$V:=(\langle\Delta\rangle\times G)\otimes_{\Bbb Z}\Bbb Q$ over $\Bbb Q$.
Then it is easily seen that  
$$\frak R:=\bigcup_{\alpha\in\Delta} (\alpha, \supp_G L_\alpha)$$
is a reduced locally extended affine root system in $V$,
extending the symmetric bilinear form having the radical $(0,G)\otimes_{\Bbb Z}\Bbb Q$.
One can now easily check that $L$ is a Lie $\frak R$-torus.
Thus we obtain:

\begin{theorem}
Any Lie $G$-torus for a torsion-free abelian group $G$
is a Lie $\frak R$-torus.
\end{theorem}

\begin{remark}
Let $G$ be a free abelian group of rank $\infty$.
Then one can construct a Lie $G$-torus by TKK construction
from a Jordan $G$-torus,
and this can be considered as
a Lie $\text A_1^{\infty}$-torus,
where
$\text A_1^{\infty}$ is an extended affine root system of type $\text A_1$
with nullity $\infty$ (in the sense of \cite{MY1}).
We note that Jordan $G$-tori 
for any torsion-free abelian group $G$
have been recently classified in \cite{AYY}.
\end{remark}

\enskip

\enskip

\section{Appendix}
\label{}

We compare the original definition of Lie tori in \cite{Y2}.
We recall $\Delta$-graded Lie algebras
introduced in \cite{BM}.
The original definition 
of a Lie torus is simply a generalization of a $\Delta$-graded Lie algebra.

\definition
Let $\Delta$ be a {\bf finite irreducible root system} and
let 
$\mathfrak g=\mathfrak h\oplus \bigoplus_{\mu\in\Delta^\red} \mathfrak g_\mu$
be a finite-dimensional split simple Lie algebra over $F$
of type $\Delta^\red$
with a split Cartan subalgebra $\mathfrak h$
and the finite irreducible reduced root system
$\Delta^\red$.

\item[(1)]
A $\Delta$-graded Lie algebra $L$ over $F$ with grading pair $(\frak g, \frak h)$ is defined as

(i) $L$ contains $\frak g$ as a subalgebra;

(ii) $L = \bigoplus_{\mu\in \Delta\cup\{0\}} L_{\mu}$, 
where $L_{\mu} = \{x\in L \mid [h,x]=\mu(h)x\ \text{for all}\ h\in\frak h\}$;

(iii) $L_0=\sum_{\mu\in\Delta}\ [L_\mu, L_{-\mu}]$.

\noindent
We also assume that 
\begin{equation}\label{stupnonred}
L_\mu\neq 0\quad\text{for all $\mu\in\Delta$}.
\end{equation}
(This is automatically true when
$\Delta$ is reduced.)

\item[(2)]
A $\Delta$-graded Lie algebra $L= \bigoplus_{\mu\in\Delta\cup\{0\}} L_{\mu}$
is called
{\it $(\Delta, G)$-graded}
if 
$L= \bigoplus_{g\in G} L^g$
is a $G$-graded Lie algebra such that 
$$\supp L:=\{g\in G\ |\ L^g\neq 0\}\ \text{generates $G$,\ \ \ \ \ and\ \ \ \ \ 
$\mathfrak g\subset L^0$. }
$$
Then we have
$$L=\bigoplus_{\mu\in\Delta\cup\{0\}}\ \bigoplus_{g\in G}\ L_{\mu}^g,$$
where $L_{\mu}^g=L_{\mu}\cap L^g$ 
since $L^g$ is an $\mathfrak h$-submodule of $L$.
Note that if $G=\{0\}$, then $L$ is just a $\Delta$-graded Lie algebra.

\item[(3)]
Let 
$Z(L)$ be the centre of $L$ and let
$\mu^{\vee}\in \mathfrak h$ for $\mu\in\Delta$
be the coroot of $\mu$.
Then $L$ 
is called a {\it division $(\Delta, G)$-graded Lie algebra}
if 
for any $\mu\in\Delta$
and any
$0\neq x\in L_{\mu}^g$,
$$\text{there exists 
$y\in L_{-\mu}^{-g}$
such that $[x,y]\equiv \mu^{\vee}$
modulo $Z(L)$.
\ \
(division property)}
$$

\item[(4)]
A division $(\Delta, G)$-graded Lie algebra
$L=\oplus_{\mu\in\Delta\cup\{0\}}\ \oplus_{g\in G}\
L_{\mu}^g$ is called a {\it Lie $G$-torus of type $\Delta$} if
$$
\text{$\dim_FL_{\mu}^g\leq 1$ for all $g\in G$ and $\mu\in\Delta$.
\ \ ($1$-dimensionality)}
$$
If $G=\Bbb Z^n$, it is called a {\it Lie $n$-torus} or 
simply a {\it Lie torus}.

\enddefinition

\remark\label{identifysub}
The assumption that $\supp L$ generates $G$ in (1) is not essential because
if $\supp L$ of a $G$-graded algebra $L$
is a proper subset of $G$,
then the subalgebra $L'$ of $L$
generated by the homogeneous spaces of degree in $G'$,
where $G'=\langle\supp L\rangle$,
is a $G'$-graded alegebra.
Moreover, $L'$ can be identified with $L$.
\endremark

\medskip

What happens if
we change the condition
$\mathfrak g\subset L^0$ in (2)
into the condition $\mathfrak h\subset L^0$\ ?

For comparison, we call the Lie $G$-torus
a {\bf normal} Lie $G$-torus,
and
the Lie $G$-torus
under the assumption $\mathfrak h\subset L^0$
a {\bf general Lie $G$-torus}.

Also, one can easily generalize both concepts
based on a finite irreducible root system $\Delta$
to the concepts based on a locally finite irreducible root system $\Delta$.
Only difference is that
$\mathfrak g=\mathfrak h\oplus \bigoplus_{\mu\in\Delta^\red} \mathfrak g_\mu$
is a locally finite split simple Lie algebra
introduced in \cite{St}
if $\Delta$ is infinite.
When we emphasize that
$\Delta$ is a locally finite irreducible root system,
we say
a general {\bf locally} Lie $G$-torus or a normal {\bf locally} Lie $G$-torus.
But we omit the term `locally' if there is no confusion,
as in Section 2.

\begin{proposition}
Two definitions in Section 2 and this section of a general (or normal) locally Lie $G$-torus
are equivalent.

\end{proposition}

\proof
Let $L$ be a general locally Lie $G$-torus defined in this section.
Then $L$ clearly satisfies (LT1-5) in Section 2.
Since
$0\neq\frak g_\mu\subset L_\mu$  for all $\mu\in\Delta^\red$,
(L6)$'$ holds
(see \eqref{stupnonred}).
If $\frak g\subset L^0$,
i.e., $L$ is normal,
then $L_\mu^0=\frak g_\mu$ for all $\mu\in\Delta^\red$, and so (L6) holds.

Next, suppose that $\mathcal L$ is a normal locally Lie $G$-torus defined in Section 2.
 From (LT3) and (LT6) we see for $\mu \in \Delta^\red$ that there
 exist  elements $e_\mu \in \mathcal L _\mu^0$, $f_\mu \in \mathcal L _{-\mu}^0$, 
 and $\mu^\vee: = [e_\mu,f_\mu]$ so that $[\mu^\vee, z] = \langle \nu,\mu \rangle z$
 for all $z \in \mathcal L _\nu^h$,
$\nu \in \Delta$ and  $h \in G$.  
Thus,  the elements
$e_\mu,f_\mu,\mu^\vee$ determine a canonical basis for
a copy of the Lie algebra $\sll_2(F)$.     
The subalgebra 
$\mathfrak g$ of $\mathcal L $ generated by the subspaces $\mathcal L _\mu^0$ for
$\mu \in \Delta^\red$  is a locally finite split simple Lie algebra with
split Cartan subalgebra 
$$\mathfrak h: = \sum_{\mu \in \Delta^\red}\
[\mathcal L _\mu^0, \mathcal L _{-\mu}^0]$$
and $\mu^\vee$ are the coroots in $\mathfrak h$.
(One can show this in the same way as the proof
of \cite[Prop.8.3]{MY1}, or see \cite[Sec.III]{St}).
Thus
$\mathcal L$ is a normal locally Lie $G$-torus defined in this section.
(Note that if $\Delta$ is finite, then $\mathfrak g$ is a finite-dimensional split simple Lie algebra.)
Suppose that $\mathcal L$ is a general locally Lie $G$-torus defined in Section 2.
By Theorem \ref{isotopicth},
$\mathcal L$ is isotopic to a normal locally Lie $G$-torus,
say $\mathcal L^{(s)}$,
and so 
$\frak h\oplus\bigoplus_{\mu\in\Delta^\red}\ \mathcal L_\mu^{s(\mu)}$
is a locally finite split simple Lie algebra,
which is a graded subalgebra of $\mathcal L=\bigoplus_{\mu\in\Delta\cup\{0\}}\ \mathcal L_\mu$
such that  $\frak h\subset \mathcal L^0$.
Thus
$\mathcal L$ is a general locally Lie $G$-torus defined in this section.
\qed

\enskip

\enskip

\enskip

\end{document}